\documentclass[twocolumn]{autart}
                \newcommand{\comments}[1]{}
\usepackage{amsmath,amsfonts,amssymb,color,enumitem}
\usepackage{graphicx}

\usepackage{caption}
\usepackage{subcaption,enumitem}
\usepackage{psfrag}
\newtheorem{theorem}{Theorem}
\newtheorem{corollary}{Corollary}
\newtheorem{proposition}{Proposition}
\newtheorem{remark}{Remark}

\newtheorem{lemma}{Lemma}

\newcommand{\rref}[1]{(\ref{#1})}
\newtheorem{assumption}{Assumption}

\definecolor{pink}{rgb}{1,0.5,0.5} % color values Red, Green, Blue

\newcommand{\mm}[1]{}
\belowdisplayskip \abovedisplayskip%
\newcommand{\ashhneg}{\renewcommand{\arraystretch}{.995}}

\newcommand{\ashalf}{\renewcommand{\arraystretch}{1.25}}
\newcommand{\ashalfplus}{\renewcommand{\arraystretch}{1.375}}
\newcommand{\asa}{\renewcommand{\arraystretch}{1.5}}
\newcommand{\asaa}{\renewcommand{\arraystretch}{1.75}}
\newcommand{\asb}{\renewcommand{\arraystretch}{2}}

\renewenvironment{thebibliography}[1]{
  \begin{oldthebibliography}{#1}
    \setlength{\itemsep}{.22em}
    \setlength{\parskip}{.22em}
}
{
  \end{oldthebibliography}
}

\begin{document}
\begin{frontmatter}

\title{ Extremum Seeking for Linear Time-Varying\\ Systems with Unknown Control Directions\thanksref{footnoteinfo}}

\vspace{-1.5em}

\thanks[footnoteinfo]{Malisoff was supported by NSF Grant   2308282 and Office of Naval Research Grant N00014-22-1-2135. Fridman was supported by  ISF-NSFC
Joint Research Program Grant/Award
Number ISF-NSFC 446/24 and \textcolor{black}{the} Chana and Heinrich
Manderman Chair at Tel Aviv University.}

%\author[First]{Jackson Knox}\hspace{1em}
\author[First]{Frederic Mazenc}\hspace{1em}
\author[Second]{Michael Malisoff}\hspace{1em}
\author[Third]{Emilia Fridman}

\address[First]{Inria Saclay, L2S-CNRS-CentraleSup\'elec, 3 rue Joliot Curie, 91192,
Gif-sur-Yvette, France (e-mail: frederic.mazenc@l2s.centralesupelec.fr)}
\address[Second]{Department of Mathematics, 303 Lockett Hall, 175 Field House Drive, Louisiana State University, Baton Rouge, LA 70803-4918, USA
(email: malisoff@lsu.edu)}

 \address[Third]{School of Electrical Engineering, Tel Aviv
University, Tel Aviv 69978, Israel (e-mail: emilia@tauex.tau.ac.il)}

\vspace{-1.25em}

\begin{abstract}
 We consider bounded extremum seeking controls for time-varying linear systems \textcolor{black}{with   uncertain coefficient matrices} and measurement uncertainty. Using a new change of variables, Lyapunov functions, and  a  comparison principle, we provide practical exponential stability bounds for the   states of the closed loop systems that hold for all nonnegative times.
 For the
first time \textcolor{black}{for linear time-varying
systems with unknown control directions, we consider bounded extremum seeking controls in the presence of \textcolor{black}{uncertain time-varying} input delays with small time-varying delay uncertainties, and we  provide reduction model controllers to compensate for the constant part of the delays}. \vspace{-.65em}
\end{abstract}

\begin{keyword}
Linear systems, estimation
\end{keyword}

\end{frontmatter}
%===============================================================================

\section{Introduction}\label{intro}
Extremum
seeking is a significant current research area, because of its ability to find
   model free, online, real time optimization methods to locate  extrema when  objective functions   include  uncertainties. Although basic  extremum seeking was studied by Leblanc in 1922 \cite{L22}, the first    stability analysis devoted to extremum seeking   appears to be in    averaging and singular perturbation methods of M.  Krstic's research team; see  \cite{S24} for a  history of  extremum seeking.  Extremum seeking is widely used, e.g.,  in aerospace systems,     networked control systems   \cite{PHTM20,SK17}, source seeking \textcolor{black}{\cite{GS23,GS24}, and  biological or chemical applications \cite{LDA19,LFD17}}. Other recent theoretical studies of extremum seeking  include
\textcolor{black}{\cite{DSEJ13,GD23,labar2022iss,MG24,SK14,suttner2022robustness}}. In particular,
 \cite{SK14,SK17}
suggested bounded extremum seeking with unknown maps in arguments of a sine or cosine, and provided bounds  for the  update rates, and \cite{SK13} provided a semiglobal analysis for extremum seeking based controllers, including   unknown control directions.

Unknown control directions   are especially challenging, since they are beyond the scope of standard control-Lyapunov functions\textcolor{black}{, and have recently been studied using} adaptive control of nonlinear systems. However, as noted in \cite{SK13}, classical parameter-adaptive solutions   exhibit undesirable transient performance, and do not achieve exponential stabilization even when no other uncertainties are present. This motivated notable   work \cite{SK13} on extremum seeking  under unknown control directions, which provided the first systematic design of extremum seeking controls for unstable plants and a new strategy for stabilization
under unknown control directions\mm{, as an alternative to
 well known Nussbaum   controllers that lack exponential stability and that lacked
 robustness to changes in the control direction}. The  main method in \cite{SK13} was Lie-brackets-based averaging.

By contrast, here we provide a new
analysis of bounded extremum seeking systems. Instead of the \textcolor{black}{Lie bracket approximation approaches from
\cite{DSEJ13,LEM22,SK14,SK17}  and the time-delay approach from \cite{ZF23,ZF24,ZF22}, we use a new delay-free change of variables, a strict Lyapunov function, and a comparison principle}, to prove practical exponential stability estimates for linear time-varying systems whose coefficient matrices may be uncertain.
\textcolor{black}{We} consider bounded extremum seeking \textcolor{black}{as introduced in \cite{SK14}} in the presence of \textcolor{black}{uncertain time-varying} delays. \textcolor{black}{This contrasts with extremum seeking results such as \cite{WGB24} whose controls are not bounded.} Robustness with respect to delays for 1D static quadratic maps was studied in \cite{YFZ23}. In the present paper, we present reduction model controllers
 that again allow us to obtain practical exponential stability estimates for our closed loop systems with bounded extremum seeking controls. This is a significant departure from   previous approaches to reduction model  control; see \cite{A82,MM16,MMN14}. Although our quadratic Lyapunov function is reminiscent of ones in \cite{MMF24a} in the absence of unknown control directions and  in the absence of drift terms, the unknown control direction problems we solve here were not amenable to   previous techniques. Also, since we do not treat the drift term   as a small perturbation term\mm{ (and so can allow unstable systems without restricting the   bounds on the coefficient  of the drift term)}, our work is beyond the scope of works such as \cite{ZF24} where the drift terms were treated as small perturbations. Finally, instead of qualitative results\mm{ for systems in closed loop with extremum seeking controls} (e.g., asymptotic stability    for small   or large enough extremum seeking parameters), our work leads to bounds on the small parameter and on the norms of the states of the closed loop systems that are valid for all nonnegative times\mm{ and which quantify the effects of the extremum seeking parameters in the upper bounding terms}.
This motivates the significantly novel approach that we present here. In summary, our main contributions are (i)
robustness with respect to measurement disturbances for fast varying \textcolor{black}{coefficient matrices}, (ii) numerical bounds on the extremum seeking parameters and ultimate bounds, and (iii) quantifying the effects of input delays in bounded extremum seeking.

We use standard notation, where $\mathbb R^n$ (resp., $\mathbb R^{p\times q}$) is the set of all vectors of real $n$-tuples (resp.,   $p\times q$  matrices),
$|\cdot|$ is the   Euclidean $2$-norm and   corresponding matrix operator norm, $|\cdot |_\infty$ is the supremum in this norm,
the dimensions of our Euclidean spaces are arbitrary, and $M_1\le M_2$ for   square matrices $M_1$ and $M_2$ of the same size \mm{in $\mathbb R^{n\times n}$ }means that  $M_2-M_1$ is positive semidefinite\textcolor{black}{, where  we recall that a symmetric matrix is called positive semidefinite (resp., definite) provided all of its eigenvalues are nonnegative (resp., positive) real numbers}. \textcolor{black}{Also},  $\sigma_{\rm min}(M)$ (resp., $\sigma_{\rm max}(M)$)  \textcolor{black}{denote} the smallest (resp., largest) singular value \cite{S98}  for    matrices $M$,
$I$ is the identity matrix, and $x_t(m)=x(t+m)$ for all $x$, $t$, and $m$ such that $t+m$ is in the domain of $x$.

\section{Class of Systems and Main Result}
\label{sec2}

We consider  single input time-varying linear systems
\begin{equation}\label{LTV}
    \dot x(t)=A(t)x(t)+B(t)u(t)\textcolor{black}{\, +\, \delta_s(t)}
\end{equation}
with   controls $u$, assuming the following on   $A$\textcolor{black}{, $B$, and $\delta_s$, where $\delta_s$ and  $\bar\delta_s$ will play key roles in Section \ref{delays2} in our study of cases with uncertain  $A$, $B$, and input delays}:
\begin{assumption}\label{as1}
The function $A: [0,+\infty)\to \mathbb R^{n\times n}$ is continuous\textcolor{black}{,    $B: [0,+\infty)\to \mathbb R^n$ is differentiable, and  $\delta_s: [0,+\infty)\to \mathbb R^n$ is piecewise continuous}. Also,
there are known constants $\bar A$, $\bar B>0$, \textcolor{black}{$\bar\delta_s$,} and  $\bar D_B$ such that
\begin{equation}\label{oldconstants}\begin{array}{l}
    |A(t)|\le \bar A,\;    |B(t)|\le \bar B,\textcolor{black}{\; |\delta_s(t)|\le \bar\delta_s,}\\ \text{and}\;  |\dot B(t)|\le \bar D_B\end{array}
\end{equation}hold for all $t\ge 0$.
\hfill$\square$\end{assumption}
\begin{assumption}\label{as2}
There are a  positive definite matrix $K\in \mathbb R^{n\times n}$,
a $C^1$  function $P : [0, + \infty) \rightarrow \mathbb{R}^{n \times n}$,
and positive constants $q$, $\underline{p}$, $\bar{p} $ and $\bar{\bar{p}}$ such that  with the choice
\begin{equation}
\label{a17v}
H(t) = A(t) -B(t) B^\top (t)K,
\end{equation}
the inequalities
\begin{equation}
\label{a19v} \begin{array}{l}
\dot{P}(t) + P(t) H(t) + H(t)^\top P(t) \leq - q P(t)
\\
\underline{p} I \leq P(t) \leq \bar{p} I,\;  |P(t)| \leq \bar{p},\; \text{and}\;
|\dot{P}(t)| \leq \bar{\bar{p}}
\end{array}\end{equation}
are satisfied   and $P(t)$ is positive definite for each $t\ge 0$.
\hfill$\square$\end{assumption}
See Remark \ref{rka} for ways in which Assumptions \ref{as1}-\ref{as2} allow cases where $A$ and  $B$ are uncertain, and Section \ref{suffic} for easily checked sufficient conditions for Assumption \ref{as2} to hold.
We fix $A$, $B$, \textcolor{black}{$\delta_s$,}   $H$, $P$, $K$, and constants satisfying   Assumptions \ref{as1}-\ref{as2} in what follows. Our   \textcolor{black}{goal} is the  practical exponential stabilization of \rref{LTV} by  a bounded extremum seeking  control, using
  the measurements
\begin{equation}\label{meas}\textstyle
y(t)=\ x^\top(t)Kx(t)+\delta(t),
\end{equation}
 where the measurement uncertainty $\delta$ satisfies:
\begin{assumption}
    \label{asdelta}
    The function $\delta: [0,+\infty)\to \mathbb R$  is piecewise continuous, and there is a known constant $\bar\delta\ge 0$ such that
     $
        |\delta(t)|\le \bar \delta$ holds for all $t\ge 0$.
\hfill$\square$\end{assumption}

\textcolor{black}{The preceding assumption is justified because in most applications, the laboratory equipment that provides measurements permits us to know bounds on the possible values of the measurement uncertainty at each time.}
Then,   for a constant value $\epsilon>0$, we   choose the bounded extremum seeking control suggested in \cite{SK14}, given by
\begin{equation}\label{uchoice}\textstyle
u(t)=\sqrt{\frac{2\pi}{\epsilon}} \cos\left(\frac{2\pi}{\epsilon} t
+ x^\top(t)Kx(t)+\delta(t)\right),
\end{equation}
and then we place conditions on $\epsilon>0$, \textcolor{black}{$\bar\delta_s$, and}  $\bar\delta$ under which the corresponding closed loop system
\begin{equation}
\label{a2v} \textstyle\begin{array}{rcl}
\dot{x}(t) &=& A(t)x(t)\textcolor{black}{\, + \, \delta_s(t)} \\&&+ B(t)
\sqrt{\frac{2\pi}{\epsilon}} \cos\left(\frac{2\pi}{\epsilon} t
+ x^\top(t)Kx(t)+\delta(t)\right)\end{array}
\end{equation}
satisfies suitable  stability conditions, when  $|x(0)|\le \sigma_0$, where $\sigma_0$ is any given positive constant.
To express our conditions and to facilitate our analysis in terms of the preceding constants, we also use
\begin{equation}
\label{MCD} \textstyle \begin{array}{rcl}
\bar{M}_C &=& \frac{1}{\sqrt{2\pi}} \left(
|HB|_\infty\! +\! |BB^\top KB|_\infty\! +\! \bar{D}_B\right)  \\&&+ \sqrt{\frac{2}{\pi^3}} \bar{A} |K|  \bar{B}^3 \epsilon+\frac{\bar B\bar \delta\sqrt{2\pi}}{\epsilon} ,\\
\lambda&=&\frac{4\pi}{4\pi-\epsilon|BB^\top K|_\infty},\;
\textcolor{black}{B^\sharp=1 + \frac{\epsilon}{4\pi}|BB^\top K|_\infty},
\end{array}
\end{equation}
 \begin{equation}\label{cis1} \!   \begin{array}{rcl}
    c_1  &=& \textcolor{black}{ 2\bar p B^\sharp  \left[ \sqrt{\epsilon} \left(\!  \bar M_C \! +\!
    \frac{2\bar B^2|K|}{\sqrt{2\pi}}\left(\bar B\bar\delta + \bar\delta_s\! \sqrt{\frac{\epsilon}{2\pi}}\right)\right) \! +\! \overline\delta_s\right]} ,\\
    c_2&=&2\bar p\left(\frac{\lambda\epsilon^2}{2\pi^2}\bar B^2\textcolor{black}{|KA+A^\top K|_\infty}|BB^\top K|_\infty\right.\\&&\left.+\,
   \textcolor{black}{ \! \lambda\bar B\sqrt{\frac{\epsilon}{2\pi}}\left(\sqrt{\frac{2\epsilon}{\pi}}|KB|_\infty\! +\! \frac{\epsilon}{\pi}|K|\bar\delta_s\right)}|BB^\top K|_\infty \right.\\&&\left.+\frac{\lambda\epsilon}{4\pi}\left[|BB^\top KH|_\infty+|HBB^\top K|_\infty\right.\right.\\&&\left.
    +|BB^\top K BB^\top K|_\infty\right]\big)+\frac{\lambda\epsilon}{\pi}\bar p\bar B\bar D_B|K|\\&&+
  \textcolor{black}{ 4\lambda\bar B|K|\left[\bar B\bar\delta+\sqrt{\frac{\epsilon}{2\pi}}\bar\delta_s\right]
  B^\sharp   \bar p}
    ,\end{array}\! \! \! \! \end{equation}
    \begin{equation}\label{cis2}\begin{array}{rcl}
c_3&=&\textcolor{black}{2\bar p
\left(\lambda^2|BB^\top K|_\infty \left[\sqrt{\frac{2\epsilon}{\pi}} |KB|_\infty+\frac{\epsilon}{\pi}|K|\bar\delta_s\right]\right.}
\\&&\textcolor{black}{\left.+2\lambda^2\sqrt{\frac{2\epsilon}{\pi}}B^\sharp  \bar{B} \bar{A} |K|
\right)} \\
c_4&=&\frac{2\lambda^3\epsilon\bar p}{\pi}
\textcolor{black}{|KA+A^\top K|_\infty} |BB^\top K|_\infty,
    \end{array}\end{equation}
\begin{equation}
\label{xis}\begin{array}{l}
\xi_s = \frac{2c_1\underline p^2}{
q\underline p^2+\sqrt{q^2\underline p^4-4c_1(c_3+2c_4)\underline p^2}},\\[1em]
\xi_l = \frac{q\underline p^2+\sqrt{q^2\underline p^4-4c_1(c_3+2c_4)\underline p^2}}{2c_3+4c_4} \; \, \text{and}\\
r_0=\frac{\xi_l-\xi_s}{\underline p^2}\left(\frac{c_3}{2}+c_4\right),\end{array}
\end{equation}
 where the positiveness of the denominator of $\lambda$ and the positiveness of $r_0$, $\xi_s$, and $\xi_l$ in \rref{xis} and the fact that $\xi_l>\xi_s$ will all follow from our next  assumption, which holds when $\epsilon>0$, \textcolor{black}{$\bar\delta_s$},   and $\bar\delta$ are small enough for any choice of our given bound $\sigma_0>0$ on the initial state $x(0)$ and for any choices of the other constants:

\begin{assumption}\label{as3}The constants defined above  are such that the four conditions
\begin{equation}\label{firstcond}\textstyle
\epsilon |BB^\top K|_\infty\le 2\pi,\vspace{-.5em}
\end{equation}
\begin{equation}
\label{secondcond}\textstyle
c_1+c_3+2c_2\le \underline pq,
\end{equation}
\begin{equation}
\label{thirdcond}
4c_1(c_3+2c_4)<\underline p^2q^2,\; \text{and}
\end{equation}
\begin{equation}
\label{fourthcond}\textstyle
2\bar p\left(\textcolor{black}{B^\sharp}\right)^2 \left(\sigma^2_0+\bar B^2\frac{\epsilon}{2\pi}\right)<\xi_l
\end{equation}
are satisfied.
\hfill$\square$\end{assumption}
In Section \ref{sec3}, we prove the following\mm{ (but see Remark \ref{weights} for other results   using different weighting terms)}:

\begin{theorem}\label{thm1}Let  Assumptions \ref{as1}-\ref{as3} hold.
Then for any constant $\xi_0$ such that
\begin{equation}\label{xi0}\begin{array}{l}
   \max\left\{\xi_s,2\bar p
   (\textcolor{black}{B^\sharp})^2\left(\sigma^2_0+\bar B^2\frac{\epsilon}{2\pi}\right)\right\}<\xi_0<\xi_l,\end{array}
\end{equation}
and for any initial state $x(0)$ for \rref{a2v} satisfying $|x(0)|\le \sigma_0$, the solution $x: [0,+\infty)\to \mathbb R^n$ of \rref{a2v} satisfies
\begin{equation}\label{conclusion}\textstyle
    |x(t)|\le \frac{\lambda}{\sqrt{\underline p}}
    \sqrt{\frac{(\xi_0 - \xi_s) \xi_le^{-r_0t } + \xi_s (\xi_l - \xi_0) }{
(\xi_0 - \xi_s)e^{-r_0t }
+
\xi_l - \xi_0
}}+  \bar B\sqrt{\frac{\epsilon}{2\pi}}
\end{equation}
for all $t\ge 0$.
\hfill$\square$\end{theorem}

\begin{remark}\label{rk21}
The motivation for  \rref{MCD}-\rref{xis} is that the $c_i$'s   in \rref{cis1}-\rref{cis2} (which are defined in terms of the required  constants from \rref{MCD}) will be used in the proof of Theorem \ref{thm1} in a differential inequality for
$V(t, R)=R^\top P(t)R$ where $P$ is from Assumption \ref{as2}, and where  $R$ will be defined using the state $x$ of the closed loop system \rref{a2v}    in the proof of the theorem. Then $\xi_s$ and $\xi_l$ will be roots of a polynomial right side of a comparison system that will be used to find an upper bound for $V(t,R(t))$. This will give \rref{conclusion} because of a relationship between $x$ and $R$.

In the special case  where $\bar A=0$, Assumption \ref{as2} implies that $|BB^\top K|_\infty\ne 0$  and $|KB|_\infty\ne 0$, because if either  supremum were zero in this case, then \rref{a19v} would give $\dot P(t)\le -q P(t)$ for all $t\ge 0$, so $P(t)\to 0$ as $t\to +\infty$, contradicting \rref{a19v}. Hence, $c_3>0$, since Assumption \ref{as1} requires that $\bar B>0$.
\hfill$\square$\end{remark}
\begin{remark}\label{rk21a}
By the positiveness of the   denominator terms in   \rref{conclusion} and the subadditivity of the square root, we can upper bound the right side of \rref{conclusion} to get
\begin{equation}\label{conclusion1}\! \! \textstyle
    |x(t)|\le \frac{\lambda}{\sqrt{\underline p}
    }
    \left(\sqrt{\frac{\xi_l(\xi_0-\xi_s)}{\xi_l-\xi_0}}e^{-0.5r_0 t}+\sqrt{\xi_s}\right)
 \! +\!  \bar B\sqrt{\frac{\epsilon}{2\pi}}
\end{equation}for all $t\ge 0$,
which gives exponential convergence
with convergence rate $0.5r_0$
to the closed ball of radius
\begin{equation}\textstyle
    B_*=
    \lambda\sqrt{\frac{\xi_s}{\underline p}}
  + \bar B\sqrt{\frac{\epsilon}{2\pi}}
\end{equation}centered at the origin,
because our formula for the rate $r_0$ in \rref{xis} is a positive constant. Therefore,  \rref{conclusion}  provides the ultimate bound
\begin{equation}\label{conc2}\textstyle
\limsup\limits_{t\to +\infty}|x(t)|\le  \lambda\sqrt{\frac{\xi_s}{\underline p}}
 + \bar B\sqrt{\frac{\epsilon}{2\pi}}.\end{equation}
Also, \textcolor{black}{when $\bar\delta_s=0$,} we can use our formulas for $\xi_s$  and $c_1$   to find a constant $\bar\zeta>0$ (not depending on $\epsilon$ or on $\bar\delta$) such that \begin{equation}\label{xisbound}\textstyle\xi_s\le\bar\zeta(1+\epsilon)\left( \sqrt{\epsilon}(1+\epsilon+\bar\delta)+\frac{\bar\delta}{\sqrt{\epsilon}}\right).\end{equation}
Therefore,   we can use\mm{ the
 subadditivity of the square root and} \rref{conc2} to find a constant $c_0>0$ (not depending on $\epsilon$ or on $\bar\delta$) such that \begin{equation}\textstyle\limsup\limits_{t\to +\infty}|x(t)|\le c_0\epsilon^{1/4}\sqrt{(1+\epsilon)\left(1+\epsilon+\bar\delta+\frac{\bar\delta}{\epsilon}\right)}.\end{equation}
For fixed values of   $t$, $\xi_s$, $\bar B$, $\epsilon$, $\underline p$, $\lambda$, $r_0$, and $\xi_l$, we can compute the partial derivative of the argument of the first square root in \rref{conclusion} with respect to $\xi_0$, to check that the right side of \rref{conclusion} is an increasing function of $\xi_0$.  Hence, choosing   $\xi_0$  closer to the left endpoint of the open interval   defined in \rref{xi0} gives a less conservative bound on $|x(t)|$ for each $t$ and set of choices of the other constants.
\mm{Different   $\xi_0$'s that satisfy \rref{xi0} give different relative sizes of the ultimate bound and the coefficient of the portion of the right side of \rref{conclusion} that exponentially decays to zero. }
\hfill$\square$\end{remark}
\begin{remark}
    \label{rka}
    Assumption \ref{as2}  is suitable for   coefficient matrices having the form
\begin{equation}\label{unkc}
    A(t)=A_0(t)+\textcolor{black}{\Delta_A}(t)\; \text{and}\; B(t)=B_0(t)+\textcolor{black}{\Delta_B}(t)
\end{equation}
for known bounded  matrix valued functions $A_0$ and $B_0$ where  the bounded   terms $\textcolor{black}{\Delta_A}$ and $\textcolor{black}{\Delta_B}$ are unknown. To see  why, notice that if $P$, $q_0>0$, and $K$ are such that
  \begin{equation}
\dot{P}(t) + P(t) H_0(t) + H_0(t)^\top P(t) \leq - q_0 P(t)\end{equation}
holds for all $t\ge 0$ with   $H_0(t)=A_0(t)-B_0(t)B^\top_0(t)K$, and if there are positive constants satisfying  the last four inequalities of \rref{a19v} for all $t\ge 0$, then \rref{a19v} will be satisfied for all $t\ge 0$ with     $q=\textcolor{black}{\omega_0q_0}$ and \rref{unkc} if \textcolor{black}{$\omega_0\in (0,1)$ is a constant such that}
\begin{equation}\! \! \! \!    \begin{array}{l}
2\bar p\left|\textcolor{black}{\Delta_A}(t)-\big[B_0(t)(\textcolor{black}{\Delta_B})^\top (t)\right.\\\left.+\textcolor{black}{\Delta_B}(t)B^\top_0(t)  \! +\!  \textcolor{black}{\Delta_B}(t)(\textcolor{black}{\Delta_B})^\top(t)\big]K\right| \le  \textcolor{black}{(1\! -\! \omega_0)q_0\underline p}\end{array} \! \! \!
\end{equation}holds
for all $t\ge 0$, which will hold under small enough bounds on $|\textcolor{black}{\Delta_A}|_\infty$ and $|\textcolor{black}{\Delta_B}|_\infty$.

\textcolor{black}{Also,  the proof of Theorem \ref{thm1} will show that the theorem remains true if the constant $\bar M_C$ and the constants $c_i$ for $i=1,2,3,4$ as defined in \rref{MCD}-\rref{cis1} are replaced by larger values and with the supremum $\epsilon |BB^\top K|_\infty$ in the denominator of the formula for $\lambda$ (as defined in \rref{MCD}) replaced by a larger value in $(0,4\pi)$ and  $D_B$ also enlarged, and  the suprema in \rref{firstcond}-\rref{conclusion} also replaced by larger values. This allows us to    generalize Theorem \ref{thm1}  to cases of the form \rref{unkc} with  additive uncertainties. This generalization is done by replacing $|HB|_\infty$ in the formula for $\bar M_C$ (which would be an unknown supremum when $A$ and $B$ are uncertain) by any upper bound for $|HB|_\infty$ (which can be computed using the known portions $A_0$ and $B_0$ of the coefficient matrices and bounds on $\textcolor{black}{\Delta_A}$ and $\textcolor{black}{\Delta_B}$) and making  similar replacements of the other supremum by larger values.
Also, since $B$ only enters our assumptions and $H$ formula through terms involving suprema or $BB^\top$, our proof also shows that we can allow $B$ to be determined only up to a sign.
We illustrate these generalizations in Section \ref{illu} below.}
\hfill$\square$\end{remark}

 \begin{remark}
The new variable $x_*=\sqrt{K}x$ transforms the closed loop system \rref{a2v} into the new system
\begin{equation}
\label{a2v1} \begin{array}{rcl}
\dot{x}_*(t) &=& A_*(t)x_*(t) \\&&+ B_*(t)
\sqrt{\frac{2\pi}{\epsilon}} \cos\left(\frac{2\pi}{\epsilon} t
+ |x_*(t)|^2+\delta(t)\right)\end{array}
\end{equation}
 where $A_*=\sqrt{K}A\sqrt{K}^{-1}$, $B_*=\sqrt{K}B$, and $\sqrt{K}$ is the principal square root. However, this transformation does not address the important question of determining which positive definite matrices $K$ in our measurements \rref{meas} allow us to satisfy the assumptions of Theorem \ref{thm1} for  given   coefficient matrices $A$ and $B$ in \rref{LTV} and thereby achieve the objectives of this paper. This motivates the sufficient conditions  in the next section.
 \end{remark}
\section{Sufficient Conditions for Assumption \ref{as2} }\label{suffic}
Before proving Theorem \ref{thm1} in Section \ref{sec3}, we explain why
Assumption \ref{as2} of the theorem holds in many cases of interest. For instance, we can prove the following, whose nondegeneracy condition \rref{PE} on $BB^\top$  can be viewed as a   persistence of excitation condition:
\begin{lemma}\label{sufficlem}
    Assume that $A: [0,+\infty)\to \mathbb R^{n\times n}$ and $B: [0,+\infty)\to \mathbb R^n$ are bounded and piecewise continuous, and that there exist constants   $\underline b>0$ and $\Delta>0$
    such that
    \begin{equation}\label{PE}\textstyle
\frac{1}{\Delta}\int_t^{t+\Delta} B(\ell) B^\top(\ell){\rm d}\ell\ge \underline b I
    \end{equation}
    for all $t\ge 0$. Choose constants $\bar A$ and $\bar B$ such that $|A(t)|\le \bar A$ and $|B(t)|\le \bar B$ hold for all $t\ge 0$, and the function \begin{equation}\textstyle
    \Gamma(t)=
    \frac{1}{\Delta} \int_t^{t+\Delta} \int_{m}^{t} B(\ell) B^\top(\ell){\rm d}\ell {\rm d}m.\end{equation} Let $\bar \Gamma\ge 0$ be any constant such that
\begin{equation}
\label{sy16}
\textcolor{black}{|\Gamma|_\infty} \leq \Delta  \bar{\Gamma}.
\end{equation}
Assume that  $K\in \mathbb R^{n\times n}$ is a positive definite matrix such that
\begin{equation}\label{sy18}  \asaa\begin{array}{l}
\textcolor{black}{|KA+A^\top K|_\infty}\left[\frac{1}{2} + \frac{2|K|^4\Delta^2\bar\Gamma^2\bar B^2}{
\underline{b}\lambda^2_{\rm min}(K)}\right]
\\+
2 \bar{A} \Delta  \bar{\Gamma} |K|^2- \frac{\underline{b}}{4}\lambda^2_{\rm min}(K) < 0.\end{array}
\end{equation}
Then, with the choices
\begin{equation}
\label{sy14a} \asaa \! \! \begin{array}{l}
P(t) = \left(1+\frac{4|K|^4\Delta^2\bar\Gamma^2\bar B^2}{\underline b\lambda^2_{\rm min}(K)}\right)\frac{K}{2}+K\Gamma(t)K,\\ \bar p=\frac{\lambda_{\rm max}(K)}{2} \left(1+\frac{4|K|^4\Delta^2\bar\Gamma^2\bar B^2}{\underline b\lambda^2_{\rm min}(K)}\right)
+\Delta\bar \Gamma|K|^2, \\
\underline p=\frac{\lambda_{\rm min}(K)}{2}  \left(1+\frac{4|K|^4\Delta^2\bar\Gamma^2\bar B^2}{\underline b\lambda^2_{\rm min}(K)}\right),\; \;  q=\frac{\underline b}{2\bar p}\lambda^2_{\rm min}(K),\\ \text{and}\; \;
\bar{\bar p}=2\bar B^2|K|^2,
\end{array}\! \!
\end{equation}
the requirements from Assumption \ref{as2} are satisfied.

\end{lemma}
{\em Proof.}
We
 consider time derivatives along solutions of the closed loop system
\begin{equation}
\label{sy1}
\begin{array}{rcl}
\dot{X}(t) & = & [A(t) - B(t) B^\top(t)K] X(t).
\end{array}
\end{equation}
 \textcolor{black}{Choose the functions}
 $
\mathcal{U}_1(X) = \frac{1}{2} X^\top KX$ and  $\mathcal{U}_2(t, X) = X^\top K\Gamma(t)K X$.
Along solutions of \rref{sy1}, we have
\begin{equation}
\label{sy5} \textstyle \ashalfplus \begin{array}{rcl}
\dot{\mathcal{U}}_1(t) &\leq& \frac{1}{2}\textcolor{black}{|KA+A^\top K|_\infty} |X(t)|^2\\&& -  X^\top(t) KB(t) B^\top(t)K X(t)
\; \, \text{and}\\
\dot{\Gamma}(t) &= &B(t) B^\top(t) - \frac{1}{\Delta}\int_t^{t+\Delta} B(\ell) B^\top(\ell){\rm d}\ell\end{array}
\end{equation}for all $t\ge 0$,   since the symmetry of $K$ gives  \begin{equation}\begin{array}{l} |X^\top(t)KA(t)X(t)|\\=\frac{1}{2}|X^\top(t)(KA(t)+A^\top (t)K)X(t)|\\\le \frac{1}{2}\textcolor{black}{|KA+A^\top K|_\infty}|X(t)|^2.\end{array}\end{equation}
Also, by  (\ref{sy5}),  and using the fact that $\lambda_{\rm min}(K^2)=\lambda^2_{\rm min}(K)$ \textcolor{black}{and \rref{sy16}}, we obtain
\begin{equation}
\label{sy7p} \! \! \! \!  \ashalfplus
\begin{array}{rcl}
\dot{\mathcal{U}}_2(t) \! & = &\!  X(t)^\top KB(t) B^\top(t) KX(t) \\\! &&\! - X(t)^\top K\frac{1}{\Delta}\int_t^{t+\Delta} B(\ell) B^\top(\ell){\rm d}\ell K X(t)
\\\! &&\! + 2 X(t)^\top K\Gamma(t) K\dot{X}(t) \\
\! & = &\!  X(t)^\top KB(t) B^\top(t) KX(t)\\\! &&\! - X(t)^\top K\frac{1}{\Delta}\int_t^{t+\Delta} B(\ell) B^\top(\ell){\rm d}\ell KX(t)
\\\! &&\! + 2 X(t)^\top K\Gamma(t) KA(t) X(t)
 \\\! &&\! - 2  X(t)^\top K\Gamma(t) KB(t) B^\top(t) KX(t)
\\
\! & \leq & \! X(t)^\top K B(t) B^\top(t) KX(t) \\\! &&\!  -  \underline{b}\lambda^2_{\scriptscriptstyle \rm min}(K) |X(t)|^2  + 2 \bar{A}|K|^2 \Delta\bar\Gamma |X(t)|^2 \\\! &&\!
+ 2 |K|^2 |X(t)| \Delta\bar\Gamma \bar{B} |B^\top(t) KX(t)|
\end{array}\! \! \! \!
\end{equation}
for all $t\ge 0$.
Moreover, we have
\begin{equation}\label{newtri}\asb\begin{array}{l}\textstyle
 2 |K|^2 |X(t)| \Delta\bar\Gamma \bar{B} |B^\top(t) KX(t)|\\=2\left\{\sqrt{\frac{\underline b \lambda^2_{\rm min}(K)}{4}}|X(t)|\right\}\\\; \; \; \; \times \left\{ \sqrt{\frac{4}{\underline b \lambda^2_{\rm min}(K)}}|K|^2\Delta\bar\Gamma \bar B|B^\top(t)KX(t)|\right\}\\
 \le  \frac{\underline b}{4}\lambda^2_{\rm min}(K)|X(t)|^2\\\; \; +\frac{4|K|^4}{\underline b\lambda^2_{\rm min}(K)}\Delta^2\bar\Gamma^2\bar B^2|B^\top(t) KX(t)|^2\end{array}\end{equation}for all $t\ge 0$, by applying the triangle inequality to the terms in \rref{newtri} in curly braces.
By using the inequality in \rref{newtri} to
  upper bound the last term in \rref{sy7p}, we obtain
\begin{equation}\label{sy10}
\begin{array}{rcl}
\dot{\mathcal{U}}_2(t)    &\le &
X(t)^\top K B(t) B^\top(t) KX(t) \\ && - \frac{3\underline{b}\lambda^2_{\rm min}(K) |X(t)|^2}{4}  + 2 \bar{A}|K|^2 \Delta\bar\Gamma |X(t)|^2 \\ &&  +  \frac{4|K|^4\Delta^2\bar\Gamma^2\bar B^2}{\underline b\lambda^2_{\rm min}(K)}|B^\top(t) KX(t)|^2
\\
  & = &   w_* X(t)^\top KB(t) B^\top(t) KX(t)
\\ && + \left(2 \bar{A} \Delta \bar\Gamma |K|^2 - \frac{3\underline{b}}{4}\lambda^2_{\rm min}(K)\right) |X(t)|^2
\end{array}
\end{equation}for all $t\ge 0$,
where
\begin{equation}
    \begin{array}{l}
    w_*=
1  +  \frac{4|K|^4\Delta^2\bar\Gamma^2\bar B^2}{\underline b\lambda^2_{\rm min}(K)}
    \end{array}
\end{equation}
We now introduce the function
$
\mathcal{U}_3(t, X) =  w_*\mathcal{U}_1(X)
+ \mathcal{U}_2(t, X)$.
Then, by    \rref{sy14a}, we have   $\mathcal U_3(t,X)=X^\top P(t)X$ \textcolor{black}{for all $t\ge 0$}.
  Hence, using (\ref{sy5}) and \rref{sy10},
 we obtain
\begin{equation}
\label{sy12a1}          \! \! \!
\begin{array}{rcl}
\dot{\mathcal{U}}_3(t)  & \leq &  w_*
\big[\frac{1}{2} \textcolor{black}{|KA+A^\top K|_\infty} |X(t)|^2 \\ && - X^\top(t) KB(t) B^\top (t)  KX(t)\big]
\\[.25em]
 & & +  w_*X^\top (t) KB(t) B^\top (t) KX(t)
\\ && + \left(2 \bar{A} \Delta \bar\Gamma |K|^2 - \frac{3\underline{b}}{4}\lambda^2_{\rm min}(K)\right) |X(t)|^2
\\[.25em]
 & = &  \left(\!
\frac{\textcolor{black}{|KA+A^\top K|_\infty}}{2} \! +\!
\frac{2|K|^4\Delta^2\bar\Gamma^2\bar B^2}{\underline b\lambda^2_{\rm min}(K)} \textcolor{black}{|KA\! +\! A^\top K\! |_\infty}\right.
\\ && \left.+ 2 \bar{A} \Delta \bar\Gamma |K|^2 - \frac{3\underline{b}}{4}\lambda^2_{\rm min}(K)
\right)
|X(t)|^2 \\
 &\le&
- \frac{\underline{b}}{2}\lambda^2_{\rm min}(K)|X(t)|^2
 \\&\leq&
- \frac{\underline{b}}{2\bar p}\lambda^2_{\rm min}(K)U_3(t,X(t))
\end{array}     \! \! \! \!
\end{equation}
for all $t\ge 0$,
where the second inequality used  \rref{sy18}. The lemma now follows from our choice of $P$ in \rref{sy14a}.\hfill$\square$

\begin{remark}\label{rka1}
Similarly to Remark \ref{rka},    \rref{PE}    allows time-varying cases \rref{unkc} with uncertainties.  To  see why, notice that   if $B_0$   admits constants $\underline b_0>0$ and $\Delta>0$ such that
\begin{equation}
\label{e78a}\textstyle
\underline{b}_0 I \leq \frac{1}{\Delta} \int_t^{t+\Delta} B_0(\ell) B_0(\ell)^\top{\rm d}\ell
\end{equation}
holds for all $t\ge 0$, then \rref{PE} holds  with $\underline b=\underline b_0/2$, provided that $|\Delta_B|_\infty$ is small enough. For instance, if \rref{e78a} holds, then we can satisfy \rref{PE} with $\underline b=\underline b_0/2$, if \begin{equation}\textstyle 2|B_0(t)||\textcolor{black}{\Delta_B}(t)|+|\textcolor{black}{\Delta_B}(t)|^2\le  \frac{\underline b_0}{2}\end{equation} holds for all $t\ge 0$. \mm{Condition \rref{e78a} can be checked in many cases, because   $B_0$ is   known; see Section \ref{illu}. Also,    \rref{sy18} can be satisfied for large enough values of $k$ and small enough values of $\Delta>0$; see our illustrations below.}
\hfill$\square$\end{remark}

\section{Proof of Theorem \ref{thm1}}
\label{sec3}
The proof has three parts. In the first part, we use
\begin{equation}\label{rz}\! \! \ashalf \begin{array}{l}
    R=\big[I\! +\! \frac{\epsilon}{4\pi}\sin\left(\frac{4\pi}{\epsilon}t
    \! +\!  2x^\top Kx\right)B(t)B^\top (t) K\big]z,\\{\rm where}\\
z(t) = x(t) -
\sqrt{\frac{\epsilon}{2\pi}} \sin\big(\frac{2\pi}{\epsilon} t
\! +\!  x^\top(t)Kx(t)\big) B(t)
\end{array}\end{equation}
to transform \rref{a2v} into a new system whose analysis is facilitated using the candidate Lyapunov function
\begin{equation}
\label{a20v}
V(t, R) = R^\top P(t) R,
\end{equation}
where $P$ is from Assumption \ref{as2}.   In the second part, we apply a comparison   argument to a differential inequality for $V$. In the third part, we use a relationship between $R$ and $x$,   to obtain the final bound \rref{conclusion}, using the fact that our choices of $\lambda$, $R$, and $z$ in \rref{MCD}  and \rref{rz} give $|z|\le \lambda |R|$.

{\em First Part.}
We use these functions:
\textcolor{black}{\begin{equation}\label{Delta1}\! \! \!  \!  \textstyle\ashalfplus\begin{array}{rcl}
\mathcal C_0(t)\!  & =  & \! \cos\left(\frac{2\pi}{\epsilon} t
+ x^\top(t) K x(t)\right),\\
\Delta_0(t)\!  &=& \! \cos\left(\frac{2\pi}{\epsilon} t
+ x^\top(t) K x(t)+\delta(t)\right)-\mathcal C_0(t),\\
\Delta_1(t)\!  &=& \! -2\mathcal C_0(t) \Delta_0(t)D(t)Kx(t)\\
\! &  & \! -2B(t)\sqrt{\frac{\epsilon}{2\pi}}\mathcal C_0(t)\delta^\top_s(t)Kx(t)+\delta_s(t)
,\\  D(t) \!  &=& \!  B(t)B(t)^\top
\end{array}\! \! \! \! \! \!
\end{equation}}
Then, our   $x$ dynamics  \rref{a2v} and   $z$ in \rref{rz} give
\begin{equation}\asa
\label{a4v1}\textstyle \! \! \! \!
\begin{array}{rcl}\asa
\dot{z}(t) \!  & = & \!  A(t)x(t) + B(t)
\sqrt{\frac{2\pi}{\epsilon}} \textcolor{black}{\mathcal C_0(t)}\textcolor{black}{\, +\, \delta_s(t)}\\\!  && \! - \dot{B}(t)
\sqrt{\frac{\epsilon}{2\pi}} \sin\left(\frac{2\pi}{\epsilon} t
+ x^\top(t) K x(t)\right)
\\
\!  & &\!   -
B(t)
\sqrt{\frac{\epsilon}{2\pi}}
\textcolor{black}{\mathcal C_0(t)}
\left(\frac{2\pi}{\epsilon}  + 2 \dot{x}^\top(t) K x(t)\right)\\\! &&\! +B(t)\sqrt{\frac{2\pi}{\epsilon}}\Delta_0(t)\\
   \! & = &  \! A(t)x(t)+B(t)\sqrt{\frac{2\pi}{\epsilon}}\Delta_0(t)\textcolor{black}{\, +\, \delta_s(t)}
\\\!  && \! - 2 B(t)\!
\sqrt{\frac{\epsilon}{2\pi}}
\textcolor{black}{\mathcal C_0(t)}
\dot{x}^\top\! (t)  K\!  x(t)
 \\\!  && \! - \sqrt{\frac{\epsilon}{2\pi}} \dot{B}(t)\sin\left(\frac{2\pi}{\epsilon} t
+ x^\top(t) K x(t)\right).
\end{array}\! \! \! \! \!
\end{equation}

Here and in the sequel, all equalities and inequalities should be understood to hold for all $t\ge 0$ along all solutions of \rref{a2v}. Hence,
   (\ref{a2v})   and \rref{Delta1} give
\begin{equation}
\label{a5v}\textstyle\! \! \! \!
\begin{array}{l}
\dot{z}(t)=B(t)\sqrt{\frac{2\pi}{\epsilon}}\Delta_0(t)+\Delta_1(t)
\\+  \left[A(t) - 2
\textcolor{black}{\mathcal C^2_0(t)}
 D(t) K
\right] x(t)
\\
 - 2 \sqrt{\frac{\epsilon}{2\pi}}  B(t)\textcolor{black}{\mathcal C_0(t)}
x^\top\! (t) A^\top\! (t) K\!  x(t)
\\- \sqrt{\frac{\epsilon}{2\pi}} \dot{B}(t)
\sin\left(\frac{2\pi}{\epsilon} t
+ x^\top(t) K x(t)\right).
\end{array}\! \! \! \!
\end{equation}

Consequently, using the half angle formula $\cos^2(l) = \frac{1}{2} [1 + \cos(2l)]$ with
$l = \frac{2\pi}{\epsilon} t + x^\top(t) K x(t)$ \textcolor{black}{to rewrite the quantity in squared brackets in \rref{a5v}}, we get
\begin{equation}
\label{a6v}\textstyle
\begin{array}{l}
\dot{z}(t)  =\sqrt{\epsilon} m_1(t, x(t))+\Delta_1(t)\\  +\left[H(t)
- \cos\left(\frac{4\pi}{\epsilon} t + 2 x^\top(t) K x(t)\right) D(t) K
\right] x(t)
\end{array}
\end{equation}
with $H(t)$ as defined by \rref{a17v} from Assumption \ref{as2},
and
\begin{equation}
\label{a7v}\textstyle\begin{array}{l}
m_1(t, x) = - \sqrt{\frac{2}{\pi}} B(t)
\textcolor{black}{\mathcal C_0}
x^\top A(t)^\top K x\\
- \frac{1}{\sqrt{2\pi}} \dot{B}(t)
\sin\left(\frac{2\pi}{\epsilon} t
+ x^\top K x\right)+B(t)
\frac{\sqrt{2\pi}}{\epsilon} \Delta_0(t).\end{array}
\end{equation}
Using   \rref{rz} and   \rref{a6v}, we therefore have
\begin{equation}
\label{b18}
\begin{array}{l}
\dot{z}(t)  = \mathcal H(t)z(t)
+ \sqrt{\epsilon} m_2(t, x(t))+\Delta_1(t)
\end{array}
\end{equation}
where
\mm{$\mathcal H(t)\! =\! H(t)
- \cos\left(\frac{4\pi}{\epsilon} t + 2 x^\top(t) K x(t)\right) D(t) K$ and }
\begin{equation}
\label{a99v}\textstyle\begin{array}{l}
\mathcal H(t)\! =\! H(t)
- \cos\left(\frac{4\pi}{\epsilon} t + 2 x^\top(t) K x(t)\right) D(t) K \\ \text{and}\; \,
m_2(t, x) =\\
\frac{1}{\sqrt{2\pi}} \mathcal H(t)
B(t)
\sin\left(\frac{2\pi}{\epsilon} t + x^\top K x\right) + m_1(t, x).\end{array}
\end{equation}
Observe for later use that our definition of $z$ in \rref{rz} gives
\begin{equation}
\label{d3}\textstyle
|x| \leq |z| + \bar{B}\sqrt{\frac{\epsilon}{2\pi}}
\end{equation}
and therefore also
\begin{equation}
\label{d2}\textstyle
|x|^2 \leq 2 |z|^2 + \bar{B}^2 \frac{\epsilon}{\pi},
\end{equation}
and that, with the choice $\bar{M}_D = 2  \bar{B} \bar{A} |K|\sqrt{2/\pi}$, we have
\begin{equation}
\label{c1a1} \asaa\begin{array}{rcl}
|m_2(t, x)| &\leq& \frac{1}{\sqrt{2\pi}} \left(
|HB|_\infty+|DKB|_\infty+\bar D_B
\right)
\\&&+ \sqrt{\frac{2}{\pi}} \bar{B} \bar{A} |K|  |x|^2+\frac{\sqrt{2\pi}}{\epsilon}\bar B\bar\delta
\end{array}\end{equation}
and so also
\begin{equation}
\label{c1a} \begin{array}{rcl}
|m_2(t, x)|
 &\le &  \bar{M}_C + \bar{M}_D |z|^2
\end{array}
\end{equation}

for all $t\ge 0$ and $x\in \mathbb R^n$,  by
using  (\ref{d2}) to upper bound the $|x|^2$  in \textcolor{black}{\rref{c1a1}} and our choice of $\bar M_C$ in \rref{MCD}.

It follows from our choice of $R$ in \rref{rz} and    \rref{b18} that
\begin{equation}\label{dotr}\ashalf
   \! \! \!  \begin{array}{rcl}
    \dot R&=&\left[I+\frac{\epsilon}{4\pi}
    \sin\left(\frac{4\pi}{\epsilon} t + 2 x^\top K x\right)D(t) K\right]\dot z\\[.28em]&&+\cos\left(\frac{4\pi}{\epsilon} t + 2 x^\top K x\right)  D(t)Kz+m_3(t,x)\\[.15em]
    &=&\left[I+\frac{\epsilon}{4\pi}
    \sin\left(\frac{4\pi}{\epsilon} t + 2 x^\top K x\right)D(t) K\right]
\mathcal H(t)z
\\[.15em]&&+\cos\left(\frac{4\pi}{\epsilon} t + 2 x^\top K x\right)  D(t)Kz+m_4(t,x)\\&=&
    \mathcal H(t)z
\! +\! \cos\left(\frac{4\pi}{\epsilon} t \! +\!  2 x^\top K x\right)   D(t)Kz
\! +\! m_5(t,x)\\[.15em]
&=& H(t)z+m_5(t,x)\\
&=& H(t)R+m_6(t,x),
    \end{array}\! \! \!
\end{equation}
where $m_3, \ldots, m_6$ are defined by:
\begin{equation}\label{mis}\ashalf\! \! \!
    \begin{array}{rcl}
    m_3(t,x)&=& \frac{\epsilon}{4\pi}\cos\left(\frac{4\pi}{\epsilon} t + 2 x^\top K x\right)4x^\top K \dot x D(t)Kz\\&&+\frac{\epsilon}{4\pi}\sin\left(\frac{4\pi}{\epsilon} t + 2 x^\top K x\right)\dot D(t)Kz,\\
    m_4(t,x)&=&
    m_3(t,x)\\&&+\left[I+\frac{\epsilon}{4\pi}\sin\left(\frac{4\pi}{\epsilon} t + 2 x^\top K x\right)D(t)K\right]\\&&\; \; \; \times \left(\sqrt{\epsilon}m_2(t,x)+\Delta_1(t)\right),\\
    m_5(t,x)&=&m_4(t,x)\\&&+\frac{\epsilon}{4\pi}
    \sin\left(\frac{4\pi}{\epsilon} t \! +\!  2 x^\top K x\right)D(t)K\mathcal H(t)z,\; \\
    m_6(t,x)&=&m_5(t,x)\\&&- \frac{\epsilon}{4\pi}H(t)\sin\left(\frac{4\pi}{\epsilon} t \! +\!  2 x^\top K x\right) \! D(t)Kz.
    \end{array}\! \! \!
\end{equation}
We next provide bounds for the functions \rref{mis}, which will make it possible to obtain a   decay estimate \textcolor{black}{for}  $V(t,R)$.
To this end, first notice that    \rref{a2v} and   \rref{d3}-\rref{d2}  give
\begin{equation}  \label{bd1a}
    \!  \! \! \!  \begin{array}{rcl}
   | m_3(t,x)|\!  &\le & \!  \big\vert\frac{\epsilon}{4\pi}\cos\left(\frac{4\pi}{\epsilon} t + 2 x^\top K x\right)4x^\top K \big\{A(t)x\\\!  && \!\left.   + B(t)\sqrt{\frac{2\pi}{\epsilon}}\cos\left(\frac{2\pi}{\epsilon} t +  x^\top K x+\delta(t)\right)\right.\\\! &&\!\left. \textcolor{black}{+\delta_s(t)}
   \right\} D(t)Kz
   \big \vert+ \frac{\epsilon}{2\pi}\bar B\bar D_B|K||z|\\
   \!  &\le & \!  \left(\frac{\epsilon \textcolor{black}{|KA+A^\top K|_\infty}}{2\pi}|x|^2\! +\! \sqrt{\frac{2\epsilon}{\pi}}|x||KB|_\infty\right.\\\! &&\! \left.\textcolor{black}{+\frac{\epsilon}{\pi}|x||K|\bar\delta_s}
   \right)\! |DK|_\infty |z|    +\frac{\epsilon}{2\pi}\bar B\bar D_B|K||z|\\\!  &\le & \!  \left[\frac{\epsilon}{2\pi}\left(2|z|^2+\bar B^2\frac{\epsilon}{\pi}\right)\textcolor{black}{|KA+A^\top K|_\infty}\right.\\ \! &&\! \left.+\left(|z|+\bar B\sqrt{\frac{\epsilon}{2\pi}}\right)
\left(\sqrt{\frac{2\epsilon}{\pi}}|KB|_\infty\right.\right.\\&&\left.\left.\textcolor{black}{+\frac{\epsilon}{\pi}|K|\bar\delta_s} \right)\right]  |DK|_\infty|z|  +\frac{\epsilon}{2\pi}\bar B\bar D_B|K||z|\\
\!  &\le & \!  \left[\frac{\epsilon}{2\pi}\left(2\lambda^2|R|^2+\bar B^2\frac{\epsilon}{\pi}\right)\textcolor{black}{|KA+A^\top K|_\infty}\right.\\ \! &&\! \left.+\left(\lambda|R|+\bar B\sqrt{\frac{\epsilon}{2\pi}}\right)
\left(\sqrt{\frac{2\epsilon}{\pi}}|KB|_\infty\right.\right.\\&&\left.\left.\textcolor{black}{+\frac{ \epsilon}{\pi}|K|\bar\delta_s} \right)\right]  |DK|_\infty\lambda|R| \\&&   +\frac{\lambda\epsilon}{2\pi}\bar B\bar D_B|K||R|.
    \end{array}\! \! \! \! \!
\end{equation}
 where  the second inequality in \rref{bd1a} used the bound \begin{equation}\begin{array}{rcl}|x^\top KA(t)x|&=&\frac{1}{2}|x^\top(KA(t)+A^\top(t)K)x|\\&\le& \frac{\textcolor{black}{|KA+A^\top K|_\infty}|x|^2}{2}\end{array}\end{equation} which follows from the symmetry of $K$, and where \rref{bd1a}  also used   \rref{firstcond} from Assumption \ref{as3} to get  $|z|\le \lambda|R|$.
 \textcolor{black}{If we now collect the coefficients of the different powers of $|R|$ from the right side of \rref{bd1a}, then we obtain}
 \begin{equation} \label{bd1}
    \!   \begin{array}{rcl}
  | m_3(t,x)|
 \! &\le  &  \! |R|\left(\frac{\lambda\epsilon}{2\pi}\bar B\bar D_B|K|\right.\\\! &&\! \left.+\frac{\lambda\epsilon^2}{2\pi^2}\bar B^2\textcolor{black}{|KA+A^\top K|_\infty} |DK|_\infty\right.\\ \! &&\left. \! +
 \textcolor{black}{\lambda\bar B\sqrt{\frac{\epsilon}{2\pi}}\left(\sqrt{\frac{2\epsilon}{\pi}}|KB|_\infty+\frac{\epsilon}{\pi}|K|\bar\delta_s\right)}\right.\\&&\left.\times |DK|_\infty\right)  + \lambda^2|DK|_\infty\left(\sqrt{\frac{2\epsilon}{\pi}}
  |KB|_\infty\right.\\&&\left.+\frac{\epsilon}{\pi}|K|\bar\delta_s\right)
  |R|^2\\\!  && \! +\frac{\epsilon}{\pi}\lambda^3\textcolor{black}{|KA+A^\top K|_\infty}|DK|_\infty|R|^3
    \end{array}\! \! \! \! \!
\end{equation}

 Similar reasoning using our \textcolor{black}{choices of $B^\sharp$ and $\Delta_1$ in \rref{MCD} and} \rref{Delta1}, our bound \rref{c1a} on $|m_2(t,x)|$, \textcolor{black}{the bound $|z|\le \lambda|R|$,}
 the relationship in \rref{a99v} between $H$ and $\mathcal H$,
 and  the bound \rref{d3} on $|x|$ gives the following bounds:
\begin{equation}\ashalf\label{bd2a}   \! \!
      \begin{array}{rcl}
   | m_4(t,x)|\! &\le &\!  |m_3(t,x)|\\\! &&\! +\left|\left[I+\frac{\epsilon}{4\pi}\sin
\left(\frac{4\pi}{\epsilon} t + 2 x^\top K x\right)D(t)K\right]\right.\\\! &&\! \left.\; \; \times\left(\sqrt{\epsilon}m_2(t,x)+\Delta_1(t)\right)\right|
\\\! &\le &\!  |m_3(t,x)|\\\! &&\!  \textcolor{black}{+B^\sharp}\left(\! \sqrt{\epsilon}\bar M_C\! +\! \sqrt{\epsilon}\bar M_D|z|^2\! +\! 2\bar B^2|K||x|\bar\delta\right.
\\\! &&\!  \left.  \textcolor{black}{+2|B|\sqrt{\frac{\epsilon}{2\pi}}\bar \delta_s|K||x|+\bar\delta_s}
\right)\\
\! &\le & \!
|m_3(t,x)|\\\! &&\!   +\textcolor{black}{B^\sharp} \left[\sqrt{\epsilon}\bar M_C + \lambda^2\sqrt{\epsilon}\bar M_D|R|^2\right.\\\! &&\! \left.  \textcolor{black}{+2\bar B|K|\left(
\bar B\bar\delta +\sqrt{\frac{\epsilon}{2\pi}}\bar\delta_s\right)}\right.\\&&\; \; \; \; \left.\textcolor{black}{\times
\left(\lambda|R| + \bar B\sqrt{\frac{\epsilon}{2\pi}}\right) + \bar\delta_s}
\right]\\
\! &=&\!  |m_3(t,x)|+\textcolor{black}{B^\sharp}\left[\textcolor{black}{\bar\delta_s}  \textcolor{black}{+\sqrt{\epsilon} (\bar M_C}\right.\\&&\! \left.\left.+
 2\bar B|K|\left(
\bar B\bar\delta \textcolor{black}{\, +\sqrt{\frac{\epsilon}{2\pi}}\bar\delta_s}\right)\frac{\bar B}{\sqrt{2\pi}}\right)\right]
\\\! &&\! \textcolor{black}{+2B^\sharp
\bar B|K|\left(
\bar B\bar\delta+\sqrt{\frac{\epsilon}{2\pi}}\bar\delta_s\right)
\lambda|R|} \\\! &&\! \textcolor{black}{+B^\sharp}\sqrt{\epsilon}\lambda^2\bar M_D|R|^2,\\
|m_5(t,x)\! &\le &\!
|m_4(t,x)|\\\! &&\! +\! \left\vert\frac{\epsilon}{4\pi}\! \sin\! \left(\frac{4\pi}{\epsilon} t \! +\!   2x^\top K x\right)D(t)K\mathcal H(t)
 z\right\vert\\
\! &\le &\!  |m_4(t,x)|\\\! &&\! +\frac{\epsilon}{4\pi}
\left(|DKH|_\infty+|DKDK|_\infty\right) \lambda|R|,\\
|m_6(t,x)|\! &\le &\!  |m_5(t,x)|\\\! &&\! +\! \left\vert \frac{\epsilon}{4\pi}\! H(t)\sin\! \left(\frac{4\pi}{\epsilon} t \! +\!   2x^\top K x\right)\! D(t)Kz\right\vert\\
\! &\le &\!  |m_5(t,x)|+\frac{\epsilon}{4\pi}|HDK|_\infty \lambda|R|.
   \end{array}\! \! \! \! \! \! \end{equation}

Moreover, Assumption \ref{as2} ensures that along all solutions of  $\dot R=H(t)R$, the function $V$ we defined in \rref{a20v} satisfies  \begin{equation}\dot V(t)\le -qV(t,R).\end{equation} Therefore, by
 collecting powers of $|R|$ on the right sides of \rref{bd1}-\rref{bd2a} and recalling our choices of the $c_i$'s   in \rref{cis1}-\rref{cis2}, it follows that along all solutions of \rref{dotr}, we have
 \begin{equation}\label{dV}\ashalf\! \! \begin{array}{rcl}
     \dot V(t)\! &\le&\!   -qV(t,R)+2R^\top(t)P(t)m_6(t,x)
     \\\! &\le&\!  -qV(t,R)\! +\! c_1|R|\! +\! c_2|R|^2\! +\! c_3|R|^3\! +\! c_4|R|^4,
 \end{array}\! \! \! \end{equation}
 where we used our choice $\bar{M}_D = 2  \bar{B} \bar{A} |K|\sqrt{2/\pi}$ to obtain the formula for $c_3$ in \rref{cis2}.
Also, by applying the triangle inequality to the terms in curly braces, we obtain
 \begin{equation}\label{t1}\begin{array}{l}
|R|=\{1\}\{|R|\} \leq \frac{1}{2} + \frac{1}{2} |R|^2
\; \, \text{and}\\
|R|^3=\{|R|\}\{|R|^2\} \leq \frac{1}{2}|R|^2 + \frac{1}{2} |R|^4.\end{array}
\end{equation} By using \rref{t1} to upper bound $|R|$ and $|R|^3$ in \rref{dV} and then collecting powers of $|R|$ and setting
\begin{equation}\label{dis}\begin{array}{l}
  d_1=
\frac{c_1}{2},\; d_2=\frac{c_1+c_3}{2}+c_2,\; \text{and}\; d_3=c_4+\frac{c_3}{2}, \end{array}\end{equation}
we obtain
    \begin{equation}\label{67}
    \dot V(t)\le -qV(t,R)+d_1+d_2|R|^2+d_3|R|^4.\end{equation}
Also, our choice of $\underline p$ from Assumption \ref{as2} and  \rref{secondcond}  give
\begin{equation}\begin{array}{l}
    \textstyle d_2|R|^2\le \frac{d_2}{\underline p}V(t,R)\le \frac{q}{2}V(t,R)\; \; \text{and}\\
    d_3|R|^4\le \frac{d_3}{\underline p^2}V^2(t,R),\end{array}
\end{equation}
which we can combine with  \rref{67} to obtain
\begin{equation}\label{dV1}\textstyle
    \dot V(t)\le -\frac{q}{2}V(t,R)+d_1+ \frac{d_3}{\underline p^2}V^2(t,R).
\end{equation}
In the next part, we use a comparison argument to find an upper bound for $V(t,R)$, which will lead to   our conclusion \rref{conclusion} of the theorem.

{\em Third Part.} In terms of the constants $d_i$ \textcolor{black}{that} we defined in \rref{dis},
consider  the comparison system
\begin{equation}
\label{n56}\textstyle
\dot{\xi} = d_1 - \frac{q}{2} \xi + \frac{d_3}{\underline p^2} \xi^2,
\end{equation}
having its initial state $\xi(0)=\xi_0$ satisfying our requirement \rref{xi0}.
By factoring the right side of \rref{n56} as a function of $\xi$  using the quadratic formula, and noting that  \rref{dis} implies that  condition \rref{thirdcond} from Assumption \ref{as3} can be rewritten as  $16 d_1d_3<  \underline{p}^2q^2$,
we can rewrite \rref{n56} as
\begin{equation}
\label{n55}\textstyle
\dot{\xi} = \frac{d_3}{\underline p^2}(\xi - \xi_s)(\xi - \xi_l)
\end{equation}
where the roots $\xi_s$ and $\xi_l$ were defined in \rref{xis}. It then follows from a partial fraction decomposition that
\begin{equation}
\label{n52}\textstyle
\frac{\xi(t) - \xi_l}{\xi(t) - \xi_s} = e^{t(d_3/\underline p^2) (\xi_l - \xi_s)} \frac{\xi(0) - \xi_l}{\xi(0) - \xi_s}
\end{equation}
and $\xi(t)\in (\xi_s,\xi_l)$ for all $t\ge 0$,
and
we can solve for $\xi(t)$ in \rref{n52} to get
\begin{equation}
\label{n51}
\begin{array}{rcl}
\xi(t) & = & \frac{
(\xi(0) - \xi_s) \xi_l + \xi_s (\xi_l - \xi(0)) e^{t (\xi_l - \xi_s) (d_3/\underline p^2)}}{
\xi(0) - \xi_s
+
(\xi_l - \xi(0))
e^{t (\xi_l - \xi_s) (d_3/\underline p^2)
}}.
\end{array}
\end{equation}

Moreover, we can use  \rref{a19v},   \rref{fourthcond},  and \rref{rz} to check that
\begin{equation}\ashalf\begin{array}{rcl}\label{74}
    V(0,R(0))&\le& \bar p|R(0)|^2\\&\le&
    \left(\textcolor{black}{B^\sharp}\right)^2 \bar p|z(0)|^2
    \\ &\le& \left(\textcolor{black}{B^\sharp}\right)^2\textcolor{black}{\bar p}\left(2|x(0)|^2+2\bar B^2\frac{\epsilon}{2\pi}\right)\\&< & \xi_l\end{array}
\end{equation}
because  $|x(0)|\le \sigma_0$\textcolor{black}{, by}  our condition \rref{xi0}  on $\xi(0)$\textcolor{black}{. Also, we can use \rref{xi0} and \rref{74} to get} $V(0,R(0))<\xi(0)$.
  By the comparison principle, we conclude that $V(t,R(t))\le \xi(t)$ for all $t\ge 0$, so  \begin{equation}\label{Rb}\textstyle|R(t)|\le \sqrt{\frac{\xi(t)}{\underline p}}\end{equation} for all $t\ge 0$,
  by our choice \rref{a20v} of $V$ and our   $\underline p$.
  Also,   \rref{firstcond} from Assumption \ref{as3} and our choice of $R$   give
  \begin{equation}\begin{array}{rcl}
      |R|&\ge& \mm{|z|-\frac{\epsilon}{4\pi}|BB^\top K|_\infty|z|\\&=&}\left(1-\frac{\epsilon}{4\pi}|BB^\top K|_\infty\right)|z|,\end{array}
  \end{equation}
    which we can combine with \rref{Rb} and our choice of $\lambda$ in  \rref{MCD} to obtain
  \begin{equation}\textstyle|z(t)|\le \frac{1}{1-\frac{\epsilon}{4\pi}|BB^\top K|_\infty}\sqrt{\frac{\xi(t)}{\underline p}}=\lambda\sqrt{\frac{\xi(t)}{\underline p}}\end{equation} for all $t\ge 0$.
    Our definition of $z$ in \rref{rz} then gives
  \begin{equation}\label{xbound}\textstyle|x(t)|\le  \lambda\sqrt{\frac{\xi(t)}{\underline p}}+\bar B\sqrt{\frac{\epsilon}{2\pi}}\end{equation} for all $t\ge 0$. The conclusion \rref{conclusion}   now follows by using \rref{n51} to upper bound the $\xi(t)$ in \rref{xbound} and by recalling our choice   $r_0=(\xi_l - \xi_s) (d_3/\underline p^2)$ from \rref{xis}.

\begin{remark}
    \label{weights}
Different   weights on the right side terms in the upper bounds in \rref{t1} produce different formulas for $d_1$, $d_2$, and $d_3$ in \rref{dis}, hence different versions of the conditions from Assumption \ref{as3} and different estimates of the ultimate bounds. To see how, notice that if $a>0$ and $b>0$ are any constants, and if
instead of \rref{t1}, we apply the triangle inequality to the terms in curly braces to get
\begin{equation}\label{t1a}\begin{array}{l}
|R|=\{\sqrt{a}\}\left\{\frac{|R|}{\sqrt{a}}\right\} \leq \frac{a}{2} + \frac{1}{2a} |R|^2
\; \, \text{and}\\[.35em]
|R|^3=\{\sqrt{b}|R|\}\left\{\frac{|R|^2}{\sqrt{b}}\right\} \leq \frac{b}{2}|R|^2 + \frac{1}{2b} |R|^4\end{array}
\end{equation}
    then we can combine terms in \rref{dV} to instead obtain \rref{67} with the new choices
$
  d_1=
\frac{ac_1}{2}$, $d_2=\frac{1}{2}\left(\frac{c_1}{a}+bc_3\right)+c_2$, \textcolor{black}{and}  $d_3=c_4+\frac{c_3}{2b}$.
Then by reasoning as in the proof of Theorem \ref{thm1}, instead of \rref{secondcond}-\rref{thirdcond}, the required conditions to obtain   \rref{conclusion} are
\begin{equation}
2d_2\le \underline pq\; \; \text{and}\; \; 16d_1d_3<q^2\underline p^2,\end{equation} with $\xi_s$ and $\xi_l$ redefined to be the roots
\begin{equation}
\textstyle
\frac{q\underline p^2\pm\sqrt{q^2\underline p^4-16d_1d_3\underline p^2}}{4d_3}.
\end{equation}
This includes  \rref{dis}, in the special case where $a=b=1$.
In Section \ref{illu}, we illustrate how choosing $a$ or $b$ differently from $1$ can reduce the ultimate bound on $|x(t)|$.
\hfill$\square$\end{remark}

\section{Using $\delta$ to Represent Measurement Delays}
\label{delays}
We can use the   $\delta(t)$ in \rref{meas}  to model \textcolor{black}{uncertain time-varying}  measurement delays\textcolor{black}{, meaning, delays in the measurements $x(t)$ that enter the output $y(t)$ in \rref{meas}}. To see how, notice that if we replace  the control $u$ in  \rref{uchoice} by
\begin{equation}\label{uchoice1}\begin{array}{l}
u_{\rm new}(t)=\\\sqrt{\frac{2\pi}{\epsilon}} \cos\left(\frac{2\pi}{\epsilon} t
+  x^\top (t-\tau(t))K x (t-\tau(t))\right),\end{array}
\end{equation}
where the piecewise \textcolor{black}{continuous} bounded function $\tau$ represents delays, then for any constant $\bar \tau>0$, the closed loop system becomes
\rref{a2v} in the special case where $\delta(t)=x^\top (t-\tau(t))Kx(t-\tau(t))-
x^\top (t)Kx(t)$, which satisfies
\begin{equation}
\begin{array}{rcl}
|\delta(t)|  &= &  \big|
x^\top (t-\tau(t))Kx(t-\tau(t))\! -\!
x^\top (t)Kx(t) \big| \\
&\le & |x^\top(t-\tau(t))K(x(t-\tau(t))-x(t))|\\&&+
 |(x^\top(t-\tau(t)-x(t))^\top Kx(t))|
 \\&\le & 2|K|B_a\left\vert\int_{t-\tau(t)}^t \dot x(\ell){\rm d}\ell\right\vert
 \\[.25em]&\le &
 2|K|B_a\left[\bar A B_a+\sqrt{\frac{2\pi}{\epsilon}}\bar B\textcolor{black}{ \, +\,  \bar\delta_s}\right]\bar\tau
\end{array}\! \! \!
\end{equation}
for all $t\in [0,\bar T]$ if $|\tau|_\infty\le \bar\tau$, where $\bar T$ is the supremum of all $t\ge 0$ such that $|x(s)|\le B_a$ for all $s\in [0,t)$ and    $B_a=c_*+B_U+\epsilon_0$ for any constant $\epsilon_0>0$, and where in terms of the notation from \textcolor{black}{\rref{conclusion1}},
\begin{equation}\begin{array}{l}
c_*=
\frac{\lambda}{\sqrt{\underline p}}
     \sqrt{\frac{\xi_l(\xi_0-\xi_s)}{\xi_l-\xi_0}}\; \text{and}\;
B_U=
\frac{\lambda}{\sqrt{\underline p}}\sqrt{\xi_s}
  +\bar B\sqrt{\frac{\epsilon}{2\pi}},
\end{array}\! \! \! \end{equation}
and  we assume that the initial function for $x$ is constant on $[-|\tau|_\infty, 0]$.
\textcolor{black}{We can assume that $\bar T>0$, by choosing $\xi_0$ close enough to $\xi_l$ so that $c_*>\sigma_0$ and using the continuity of $x(s)$.} By the definition of the supremum and the continuity of $x$, it follows that  if $\bar  T$ is finite, then $|x(\bar  T)|=B_a$. We now argue by contradiction.  If the assumptions of Theorem \ref{thm1} hold with \begin{equation}\label{bardel}\textstyle\bar \delta=2|K|B_a\left[\bar A B_a+\sqrt{\frac{2\pi}{\epsilon}}\bar B\textcolor{black}{ \, +\,  \bar\delta_s}\right]\bar\tau, \end{equation}
then
we can apply the   argument  in the proof of Theorem \ref{thm1} except only for  $t\in[0,\bar  T)$    to obtain \rref{conclusion1} for all $t\in [0,\bar  T)$. Hence,  if $\bar  T<+\infty$, then  \rref{conclusion1} would give $|x(\bar  T)|<B_a$, since $r_0>0$. This contradiction implies that $\bar  T=+\infty$.

Therefore, if   the assumptions of Theorem \ref{thm1} hold with   \rref{bardel}, then the conclusion of the theorem remains true when the original undelayed control $u$ from \rref{uchoice} is replaced by   \rref{uchoice1}.  Assumptions \ref{as1}-\ref{as3} are satisfied with  \rref{bardel} and small enough $\epsilon>0$ and $\bar\delta>0$, provided
\begin{equation}\label{taubound}
    \begin{array}{l}
    |\tau|_\infty\le \frac{\bar\delta}{2|K|B_a\left[\bar A B_a+\sqrt{2\pi/\epsilon}\bar B\textcolor{black}{ \, +\,  \bar\delta_s}\right]}
    \end{array}
\end{equation}
 (which is similar to the constraint $\tau=O(\sqrt{\epsilon})$ in \cite{ZF23,ZF24}).  Also, by  reasoning as in \cite[Section IV]{MMF24a} (which was confined to  cases where $A=0$ and $B=I$), it follows that a restriction on the bound for  $\tau$ cannot be removed.

\section{Reduction Model Approach}\label{delays2}

  \textcolor{black}{We can also use} a reduction model approach to prove a  result \textcolor{black}{for uncertain time-varying delays} \textcolor{black}{that} incorporates the effects of the delays into an upper  bound on the norm $|x(t)|$ of the state   that holds for all nonnegative times, as follows. We study the system
\begin{equation}\label{delayeq}
\dot x(t)=A(t)x(t)+B(t)u(t-\tau\textcolor{black}{(t)})
\end{equation}with \textcolor{black}{$A$ and $B$ having the form \rref{unkc} (with $A_0$ and $B_0$ being known matrix valued functions, but $\textcolor{black}{\Delta_A}$ and $\textcolor{black}{\Delta_B}$ being unknown bounded matrix valued      functions) and an unknown time-varying delay $\tau(t)$ having the form\begin{equation}\label{tdecomp}
    \tau(t)=\hat\tau+\Delta_\tau(t)
\end{equation}
where $\hat\tau> 0$ is a known constant  and $\Delta_\tau$ is an unknown piecewise continuous bounded function}, with the control
 \begin{equation}
\label{uchoice2a}\begin{array}{l}
u(t) =
\sqrt{\frac{2\pi}{\epsilon}} \cos\left(\frac{2\pi}{\epsilon} t
+  \zeta^\top (t)K\zeta(t)\right),
\\
{\rm where}\; \zeta(t) = x(t) +  \mathcal T(t,\zeta_t)\; \text{and}\\ \mathcal{T}(t, \zeta_t) = \textcolor{black}{\int_{t-\hat\tau}^{t} \Phi_{E}(t-\hat\tau, m) B_0(m + \hat\tau) u(m){\rm d}m}\\ \text{for all}\;  t\ge 0
\;  \text{and}\; \zeta(t)=\zeta(0)\; \text{for all} \; t<0
\end{array}
\end{equation}
  which is a generalization of the reduction model transformation from \cite{MM16} (which did not cover extremum seeking)
 where  $\Phi_E$ is the state transition matrix for $E(t)=A_0(t+\textcolor{black}{\hat \tau})$.
\textcolor{black}{We} now set
\begin{equation}\label{94}B_{\rm new}(t) = \Phi_{E}(t-\textcolor{black}{\hat \tau}, t) B_0(t + \textcolor{black}{\hat \tau}),\end{equation}
   \begin{equation}  \label{Kstar}\! \! \! \begin{array}{l}K_*=\sup\{|\Phi_E(t-\textcolor{black}{\hat \tau},m)|
: t - \textcolor{black}{\hat\tau}\le  m \le  t, t\ge 0\},\end{array} \end{equation} \textcolor{black}{and}\textcolor{black}{
\begin{equation}\label{barMtau}\! \! \! \begin{array}{l}
      \bar {\mathcal M}_\Delta=\left(\bar\Delta_A \hat\tau K_* |B_0|_\infty+\bar\Delta_B \right)\sqrt{\frac{2\pi}{\epsilon}}\\
    +\, 2\bar B\sqrt{\frac{2\pi}{\epsilon}}M_\tau |K|\big[ |A|_\infty M_\tau+\left(|B_{\rm new}|_\infty+2\bar B\right.\\\left.\left.+\, \bar\Delta_A \hat\tau K_* |B_0|_\infty +\bar\Delta_B\right)\sqrt{\frac{2\pi}{\epsilon}}\, \right]|\Delta_\tau|_\infty
    ,\end{array}\end{equation}where $\overline B$, $\bar\Delta_A$, $\bar\Delta_B$, and $\bar\Delta_\tau$ are any known upper bounds for $|B|_\infty$, $|\Delta_A|_\infty$, $|\Delta_B|_\infty$, and $|\Delta_\tau|_\infty$ respectively,
    \begin{equation}\textstyle
    M_\tau=\frac{\lambda}{\sqrt{\underline p}}\sqrt{\xi_0}+|B_{\rm new}|_\infty\sqrt{\frac{\epsilon}{2\pi}},
\end{equation}
 $\xi_0>0$ will be the initial state of a comparison system, and $\underline p$, $\lambda$, and $\epsilon$ will satisfy our assumptions above. These constants will all be bounded, because of our assumptions from our theorem of this section.}

\textcolor{black}{Our reduction model theorem is as follows,} where $r_0>0$ is the exponential convergence rate from \rref{xis} which is computed with $B$ in Theorem \ref{thm1} replaced by $B_{\rm new}$ and therefore with $\bar B$ and $\bar D_B$ replaced by $|B_{\rm new}|_\infty$ and $|\dot B_{\rm new}|_\infty$, respectively\textcolor{black}{, and where   \rref{newepscond} can be interpreted as a smallness condition on $\epsilon$ or $\textcolor{black}{\Delta_A}$ that implies that we cannot compensate for an arbitrary   $\hat\tau>0$}\mm{, and where
$|\Delta_\tau |_\infty=O(\epsilon)$ follows from \rref{barMtau} and our use of Assumption \ref{as3} in Theorem \ref{prop1} with $\bar {\mathcal M}_\Delta=\bar \delta_s$}
(but see Proposition \ref{as2remark} below, for sufficient conditions  for the assumptions of Theorem \ref{prop1} to hold):
\begin{theorem}\label{prop1}
   \textcolor{black}{ Let  $\sigma_0>0$ be  a constant,   $A:\mathbb R\to \mathbb R^{n\times n}$ and $B:[0,+\infty)\to \mathbb R^n$ be bounded continuous functions of the form \rref{unkc}, with $B$ being differentiable with a bounded derivative and $\textcolor{black}{\Delta_A}$ and $\textcolor{black}{\Delta_B}$ being bounded and continuous, and $\tau:[0,+\infty)\to [0,+\infty)$ have the form \rref{tdecomp} as above.}
    Let Assumptions \ref{as1}-\ref{as3} hold with the preceding    $A(t)$  but with $B(t)$ in Assumptions \ref{as1}-\ref{as3} replaced by $B_{\rm new}(t)$ and with $\bar\delta=0$ and \textcolor{black}{$\bar\delta_s=\bar{\mathcal M}_\Delta$}, and  let \textcolor{black}{$\epsilon\in (0,1)$, $\xi_s$, $\xi_l$,} and  $K$ satisfy these assumptions with this replacement.
    Let
    \begin{equation}\textstyle
        \mathcal L=2\bar p\left(1\! +\! \frac{\epsilon}{4\pi}|B_{\rm new}B^\top_{\rm new} K|_\infty\right)^2\left(\sigma^2_0\! +
        \! \frac{\epsilon|B_{\rm new}|_\infty^2}{2\pi}\right).
    \end{equation}
 Then, for  each initial state $x(0)$ for \rref{delayeq} such that  $|\zeta(0)|\le \sigma_0$ and  each $\xi_0 \in (\max\{\xi_s,\mathcal L\},\xi_l)$, where $\xi_s$ and $\xi_l$ are computed using   \rref{xis} except with $B$ replaced by $B_{\rm new}$,
    the solution $x(t)$ of
\rref{delayeq} in closed loop with the control   \rref{uchoice2a}
    satisfies
    \begin{equation}
\label{kn10}
\begin{array}{rcl}
|x(t)|
& \le & \frac{\lambda}{\sqrt{\underline p}} \sqrt{\frac{
(\xi_0 - \xi_s) \xi_le^{-r_0t } + \xi_s (\xi_l - \xi_0) }{
(\xi_0 - \xi_s)e^{-r_0t }
+
\xi_l - \xi_0
}}\\&&+\,   |B_{\rm new}|_\infty\sqrt{\frac{\epsilon}{2\pi}}
+  |\mathcal T(t,\zeta_t)|
\end{array}
\end{equation}for all $t\ge 0$.
Also, \textcolor{black}{if the condition
\begin{equation}\! \! \!
\label{newepscond}\begin{array}{l}\textcolor{black}{
   \frac{2|K|\hat\tau\sqrt{2\epsilon}}{\sqrt{\pi}}
      K_* \bar\Delta_A |B_0|_\infty\left(\lambda\sqrt{\frac{\xi_s}{\underline p}}\! +\! |B_{\rm new}|_\infty\sqrt{\frac{\epsilon}{2\pi}}\right)<1}\end{array} \!
\end{equation}
is also satisfied, then}
we can  find positive constants $c_\mathcal T$
and $d_{\mathcal T}>0$ \textcolor{black}{(both not depending on $\epsilon$ or $\bar\delta_s$ or $\textcolor{black}{\Delta_B}$)}
such that  \textcolor{black}{with the choice
\begin{equation}
\label{epsilonsharp}\epsilon^\sharp=\left(1+\bar\Delta_B\right)\left(\epsilon^{1/4}+\sqrt{\bar\delta_s}\right)+\sqrt{\epsilon}\bar\delta_s,
\end{equation}}
the inequalities
\begin{equation}\label{nucl}\! \! \begin{array}{l}
\limsup\limits_{t\to +\infty}|\mathcal T(t,\zeta_t)|\le  c_{\mathcal T}\textcolor{black}{\epsilon^\sharp}\; \, \text{and}\\ \limsup\limits_{t\to +\infty}|x(t)|\le d_{\mathcal T}\textcolor{black}{\epsilon^\sharp}
\end{array}\! \! \! \end{equation} are satisfied.
\hfill$\square$\end{theorem}

{\em Proof.} The proof entails applying the reduction model approach to \rref{delayeq}.
\textcolor{black}{Consider any solution $x(t)$ of the closed loop system from Theorem \ref{prop1} such that $|\zeta(0)|\le \sigma_0$. Then
\rref{uchoice2a}
gives}
\begin{equation}
\label{kn2}\! \! \!   \textcolor{black}{
\begin{array}{rcl}
\dot{\zeta}(t)  &=&  A(t) x(t) + B(t) u(t - \tau(t))
\\&&+ E(t - \hat \tau)
\mathcal T(t,\zeta_t)  - B_0 (t) u(t - \hat \tau )
\\
&& + \Phi_{E}(t-\hat \tau, t) B_0(t + \hat \tau) u(t),
\end{array}}\! \!
\end{equation}\textcolor{black}{\noindent where throughout the proof, all equalities and inequalities should be understood to hold for all $t> 0$ unless otherwise indicated.}
Since \begin{equation}E(t-\textcolor{black}{\hat \tau})=A\textcolor{black}{_0}(t),\end{equation} we can therefore use   \rref{uchoice2a}-\rref{94} to get
 \begin{equation}
\label{kn3}\! \! \!
\textcolor{black}{\begin{array}{rcl}
\dot{\zeta}(t) & = & A(t) \zeta(t) + B_{\rm new}(t) u(t) \textcolor{black}{\, +\, \mathcal M(t)}
\end{array}}\! \! \!
\end{equation}
  where
\begin{equation}\! \! \! \label{calM}\textcolor{black}{\begin{array}{rcl}
\mathcal M(t)&=&B(t)[u(t-\tau(t))-u(t-\hat \tau)]
\\&&-\textcolor{black}{\Delta_A}(t)\mathcal T(t,\zeta_t)+\textcolor{black}{\Delta_B}(t)u(t-\hat \tau).\end{array}}\end{equation} \textcolor{black}{Also,
our lower bound $\max\{\xi_s, \mathcal L\}$
on $\xi_0$  and the facts that  $\lambda>1$  and $\bar p\ge \underline p$ give $|\zeta(0)|\le \sigma_0<M_\tau$. Therefore, by the continuity of $\zeta$, we can find a $t_*>0$  such that $|\zeta|_{\scriptscriptstyle [0,t_*]}<M_\tau$, where   $|\cdot |_J$   indicates the supremum over any interval $J$. Therefore, if we define $\bar T$ by $\bar T=\sup\{t\ge 0: |\zeta|_{[0,t]}<M_\tau\}$, then $\bar T$ is positive or is $+\infty$.
Also, using the  function \begin{equation}\label{zetak}\zeta_K(m)=\zeta^\top(m)K\zeta(m),\end{equation} it follows that for all $\ell\in [0,\bar T)$, we can use \rref{kn3} to get
\begin{equation}\label{Mstar}\asa
    \begin{array}{l}
    |B(\ell)[u(\ell-\tau(\ell))-u(\ell-\hat \tau)]|\\
    \le \bar B\sqrt{\frac{2\pi}{\epsilon}}|\zeta_K(\ell-\tau(\ell))-\zeta_K(\ell-\hat \tau)|\\
    \le 2\bar B\sqrt{\frac{2\pi}{\epsilon}}M_\tau |K|\int_{(\ell-\max\{\tau(\ell),\hat\tau\})^+}^{(\ell-\min\{\tau(\ell),\hat\tau\})^+}\left[ |A|_\infty|\zeta(s)|\right.\\\left.\; \; \; +\left(|B_{\rm new}|_\infty+2|B|_\infty+|\textcolor{black}{\Delta_A}|_\infty \hat \tau K_* |B_0|_\infty\right.\right.\\\left.\left.\; \; \; +|\textcolor{black}{\Delta_B}|_\infty\right)\sqrt{2\pi/\epsilon}\right]{\rm d}s\\
    \le 2\bar B\sqrt{\frac{2\pi}{\epsilon}}M_\tau |K|\left[ |A|_\infty M_\tau+\left(|B_{\rm new}|_\infty+2|B|_\infty\right.\right.\\\left.\left.+|\textcolor{black}{\Delta_A}|_\infty \hat\tau K_* |B_0|_\infty +|\textcolor{black}{\Delta_B}|_\infty\right)\sqrt{\frac{2\pi}{\epsilon}}\, \right]|\Delta_\tau|_\infty,
    \end{array}\! \! \! \!
\end{equation}because $\cos$ has the global Lipschitz constant $1$, where  $r^+=\max\{r,0\}$ for each $r\in \mathbb R$. Therefore, our bound $2\pi/\epsilon$ on   $u$ gives $|\mathcal M(t)|\le \bar{\mathcal M}_\Delta$ for all $t\in [0,\bar T)$.}

\textcolor{black}{
Next, notice that
\rref{kn3} has the form \rref{a2v}  with $B$ replaced by $B_{\rm new}$ \textcolor{black}{and $\delta_s(t)=\mathcal M(t)$, and its control from \rref{uchoice2a} agrees with the one from Theorem \ref{thm1}, when $\bar\delta=0$}.
Therefore, applying the argument from the proof of Theorem \ref{thm1} except only for   $t\in [0,\bar T)$, with  $\delta_s=\mathcal M$,  $\bar\delta=0$, and   $B$ replaced by $B_{\rm new}$,  implies that for all $t\in [0,\bar T)$, we have \begin{equation}\label{ins}\begin{array}{l}
|\zeta(t)| \le\\ \frac{\lambda}{\sqrt{\underline p}} \sqrt{\frac{
(\xi_0 - \xi_s) \xi_le^{-r_0t } + \xi_s (\xi_l - \xi_0) }{
(\xi_0 - \xi_s)e^{-r_0t }
+
\xi_l - \xi_0
}}+|B_{\rm new}|_\infty\sqrt{\frac{\epsilon}{2\pi}}.
\end{array}\end{equation}
 Also, if $\bar T<+\infty$, then  $|\zeta(\bar T)|<M_\tau$, since $\zeta$ is continuous and   the first right side term of \rref{ins}  is strictly decreasing in $t$   and $\bar T>0$. By again using the continuity of $\zeta$, we could then find a constant $\epsilon_0>0$ such that $|\zeta(\bar T+\ell)|<M_\tau$ for all $\ell\in [0,\epsilon_0]$, which contradicts our choice   of $\bar T$. Hence, $\bar T=+\infty$, so \rref{ins} holds for all $t\ge 0$.
Since $x(t)=\zeta(t)-\mathcal T(t,\zeta_t)$, it follows from   the
 triangle inequality that
\rref{kn10} holds for all $t\ge 0$.
  \mm{This proves the first statement of the theorem.}}

 \textcolor{black}{We next derive the required constant   $c_{\mathcal T}$, in \rref{nucl}, and then we use $c_{\mathcal T}$ to find $d_{\mathcal T}$.
 First note} that we can apply the sum rule for cosine    in our formula for $\mathcal T $, to  obtain
 \begin{equation}\label{decomp}\! \! \! \begin{array}{l}
\mathcal{T}(t, \zeta_t) = \sqrt{\frac{2\pi}{\epsilon}}\mathcal{T}_1(t, \zeta_t) - \sqrt{\frac{2\pi}{\epsilon}}\mathcal{T}_2(t, \zeta_t),\end{array}\end{equation}\textcolor{black}{for all $t\ge \hat\tau$,}
where \textcolor{black}{in terms of our function \rref{zetak},}
\begin{equation}\label{decomps}\! \! \! \begin{array}{l}
\mathcal{T}_1(t, \zeta_t)\\= \int_{t-\textcolor{black}{\hat\tau}}^{t}
\cos\left(\frac{2\pi}{\epsilon} m\right)
\lambda(t, m)
\cos\left(\textcolor{black}{\zeta_K(m)}\right)
{\rm d}m\textcolor{black}{,}\\
\mathcal{T}_2(t, \zeta_t)\\ =
\int_{t-\textcolor{black}{\hat\tau}}^{t} \sin\left(\frac{2\pi}{\epsilon} m\right)\lambda(t, m)
\sin\left(\textcolor{black}{\zeta_K(m)}\right){\rm d}m,\\
 \text{and}\;
 \lambda(t, m)= \textcolor{black}{\Phi_{E}(t-\hat\tau, m) B_0(m + \hat\tau)}.\end{array}\end{equation}
If we now integrate $\mathcal T_1$ by parts and use  \rref{kn3} and our choice \rref{uchoice2a} of   $u$,
then we get
\begin{equation}
\label{kn15}    \! \!   \! \!   \asa
\begin{array}{l}
\mathcal{T}_1(t, \zeta_t)     =     \frac{\epsilon}{2\pi}
\left[\sin\left(\frac{2\pi t}{\epsilon}\right) \lambda(t,t)
\cos\left(\textcolor{black}{\zeta_K(t)}\right)
 \right.\\\left.- \sin\left(\frac{2\pi (t  -  \textcolor{black}{\hat\tau})}{\epsilon}\right)\lambda(t,t \! -\!  \textcolor{black}{\hat\tau})
\cos\left(\textcolor{black}{\zeta_K(t-\hat\tau)}
 \right)\right]
\\
      - \frac{\epsilon}{2\pi} \! \int_{t-\textcolor{black}{\hat\tau}}^t
\! \sin\left(\frac{2\pi m}{\epsilon}\right)
\! \frac{\partial \lambda}{\partial m}(t, m) \cos(\textcolor{black}{\zeta_K(m)})
{\rm d}m
\\
       + \frac{ \epsilon}{\pi} \int_{t-\textcolor{black}{\hat\tau}}^t
\sin\left(\frac{2\pi m}{\epsilon}\right)
\lambda(t, m) \sin(\textcolor{black}{\zeta_K(m)})\\\; \; \; \; \; \; \times \zeta(m)^\top K\dot{\zeta}(m)
{\rm d}m
\end{array}\! \! \! \! \!   \end{equation}
and therefore also
\begin{equation}
\label{kn15a}         \! \! \! \! \asaa
\begin{array}{l}
\mathcal{T}_1(t, \zeta_t)
    =     \frac{\epsilon}{2\pi}
\left[
\sin\left(\frac{2\pi t}{\epsilon}\right) \lambda(t,t) \cos\left(\textcolor{black}{\zeta_K}(t)\right)
\right.\\\left.- \sin\left(\frac{2\pi (t  -  \textcolor{black}{\hat\tau})}{\epsilon}\right) \lambda(t,t  -  \textcolor{black}{\hat\tau})
\cos\left( \textcolor{black}{\zeta_K(t - \hat\tau))}\right)
\right]
\\
    - \frac{\epsilon}{2\pi}\!  \int_{t-\textcolor{black}{\hat\tau}}^t
\sin\left(\frac{2\pi m}{\epsilon}\right)
\frac{\partial \lambda}{\partial m}(t, m) \cos(\textcolor{black}{\zeta_K(m)})
{\rm d}m
\\
       + \frac{ \epsilon}{\pi} \int_{t-\textcolor{black}{\hat\tau}}^t
\sin\left(\frac{2\pi m}{\epsilon}\right)
\lambda(t, m) \sin(\textcolor{black}{\zeta_K(m)}) \\\; \; \; \; \; \; \times \zeta(m)^\top K
A(m) \zeta(m)
{\rm d}m
\\    +  \textcolor{black}{\frac{\epsilon}{\pi}}  \int_{t-\textcolor{black}{\hat\tau}}^t
\sin\left(\frac{2\pi m}{\epsilon}\right)
\lambda(t, m) \sin(\textcolor{black}{\zeta_K(m)}) \\\; \; \;  \times \zeta^\top(m) K
\left[\textcolor{black}{\sqrt{\frac{2\epsilon}{\pi}}}B_{\rm new}(m)
\cos\left(\frac{2\pi}{\epsilon} m
\!  +\!   \textcolor{black}{\zeta_K(m)}\right) \right.\\\left.\; \; \; \textcolor{black}{+\mathcal M(m)}\right]
{\rm d}m.
\end{array}    \!  \! \! \! \!    \! \! \! \! \!
\end{equation}Therefore,
\begin{equation}
\label{kn17}
\begin{array}{l}
|\mathcal{T}_1(t, \zeta_t)|
 \leq  \frac{\epsilon}{2\pi} \left(
|\lambda(t,t)| + |\lambda(t,t - \textcolor{black}{\hat\tau})|\right)
\\+\,  \frac{\epsilon}{2\pi} \int_{t-\textcolor{black}{\hat\tau}}^t
\left|\frac{\partial \lambda}{\partial m}(t, m)\right|{\rm d}m
\\
 + \frac{|K| \bar{A}\epsilon}{\pi} \int_{t-\textcolor{black}{\hat\tau}}^t
|\lambda(t, m)||\zeta(m)|^2
{\rm d}m
\\+\,  |K| \textcolor{black}{\mathcal B_\tau(t)}
\textcolor{black}{\frac{\epsilon}{\pi}}
\int_{t-\textcolor{black}{\hat\tau}}^t|\lambda(t, m)| |\zeta(m)|
{\rm d}m
\end{array}
\end{equation}
\textcolor{black}{holds, where $\mathcal B_\tau(t)=\sqrt{2\epsilon/\pi}|B_{\rm new}|_\infty+|\mathcal M|_{[\max\{0,t-\hat\tau\},t]}$.
Hence,  we can use  \cite[Lemma C.4.1]{S98}   and
  (\ref{ins}) to get}
\[\ \! \! \! \! \asa
\begin{array}{l}
|\mathcal{T}_1(t, \zeta_t)|
 \leq  \frac{\epsilon}{2\pi} C_*\\[.05em]
+ \frac{|K|\bar A\epsilon\textcolor{black}{ |B_0|_\infty\hat\tau}}{\pi} K_*
\left(
\frac{\lambda}{\sqrt{\underline p}} \sqrt{\frac{
(\xi_0 - \xi_s) \xi_l\textcolor{black}{r_*(t)} + \xi_s (\xi_l - \xi_0) }{
(\xi_0 - \xi_s)\textcolor{black}{r_*(t)}
+
\xi_l - \xi_0
}}\right.\\+|B_{\rm new}|_\infty\sqrt{\frac{\epsilon}{2\pi}}
\bigg)^2
+\textcolor{black}{\mathcal B_\tau(t)}\textcolor{black}{|B_0|_\infty|K|\frac{\hat\tau\epsilon}{\pi}K_*}\!
\\\times \left(\!
\frac{\lambda}{\sqrt{\underline p}} \sqrt{\frac{
(\xi_0 - \xi_s) \xi_l\textcolor{black}{r_*(t)} + \xi_s (\xi_l - \xi_0) }{
(\xi_0 - \xi_s)\textcolor{black}{r_*(t)}
+
\xi_l - \xi_0
}}+|B_{\rm new}|_\infty\sqrt{\frac{\epsilon}{2\pi}}
\right)
\end{array} \! \! \! \!
\]  \textcolor{black}{for all $t\ge \hat \tau$,} where \textcolor{black}{$r_*(t)=e^{-r_0(t-\hat\tau)}$ and}
\[C_*=\textcolor{black}{|B_0|_\infty} \left(K_*\! +\! 1\right)\! +\! \textcolor{black}{\hat\tau}\bar A\textcolor{black}{|B_0|_\infty} K_*\! +\! \textcolor{black}{\hat\tau}|\dot B\textcolor{black}{_0}|_\infty K_*,\] by  using the fact that the argument of the first square root in \rref{kn10} is an decreasing function of $t$.
Similar reasoning with the roles of   $\sin$ and $\cos$ interchanged gives
\begin{equation}
\label{kn18a}\! \! \!
\begin{array}{l}
\limsup\limits_{t\to +\infty}|\mathcal{T}_i(t, \zeta_t)|
   \leq   \frac{\epsilon}{2\pi} C_*
\\+ \frac{|K|\bar A\epsilon\textcolor{black}{|B_0|_\infty\hat\tau}}{\pi} K_*
 \left(
\frac{\lambda}{\sqrt{\underline p}}\sqrt{\xi_s}  +| B_{\rm new}|_\infty\sqrt{\frac{\epsilon}{2\pi}}
\right)^2
\\
+ |K|\textcolor{black}{
\frac{\hat\tau\epsilon}{\pi}}
K_*\textcolor{black}{\mathcal B_\tau(t)}|B\textcolor{black}{_0}|_\infty \left(
\frac{\lambda\sqrt{\xi_s}}{\sqrt{\underline p}} \! +\! | B_{\rm new}|_\infty\sqrt{\frac{\epsilon}{2\pi}}
\right)
\end{array} \! \! \! \! \! \! \! \! \!
\end{equation}for $i=1,2$.
Since  \textcolor{black}{\rref{cis1} and \rref{xis} yield} a constant $\bar\xi_s$ \textcolor{black}{that does not depend on $\epsilon$ or $\bar\delta_s$}  such that \begin{equation}\label{xisbound1}\textcolor{black}{\xi_s\le \bar\xi_s (\sqrt{\epsilon}+\bar\delta_s)}\end{equation} when $\bar\delta=0$ and $\epsilon\in (0,1)$, \textcolor{black}{we can then combine the decomposition from  \rref{decomp} with our choice \rref{calM} of $\mathcal M$ and \rref{kn18a} to obtain}
\begin{equation}\label{fTb}
    \textcolor{black}{\begin{array}{l}
    \limsup\limits_{t\to +\infty}|\mathcal T(t,\zeta_t)|\\\le 2\sqrt{\frac{2\pi}{\epsilon}}\left[\frac{\epsilon}{2\pi}C_*+\frac{|K|\bar A\epsilon\hat\tau|B_0|_\infty}{\pi}  K_*\left(\frac{\lambda}{\sqrt{\underline p}}\sqrt{\xi_s}\right.\right.\\+|B_{\rm new}|_\infty\sqrt{\frac{\epsilon}{2\pi}}\bigg)^2+
    |K|\frac{\hat\tau\epsilon}{\pi}K_*\bigg(\sqrt{\frac{2\epsilon}{\pi}}|B_{\rm new}|_\infty\\+2|B|_\infty\sqrt{\frac{2\pi}{\epsilon}}+|\textcolor{black}{\Delta_A}|_\infty \limsup\limits_{t\to +\infty}|\mathcal T(t,\zeta_t)|\\\left.\left.+|\textcolor{black}{\Delta_B}|_\infty\sqrt{\frac{2\pi}{\epsilon}}\right)|B_0|_\infty
\left(\lambda\sqrt{\frac{\xi_s}{\underline p}}+|B_{\rm new}|_\infty \sqrt{\frac{\epsilon}{2\pi}}\right)\right]
    \end{array}}\! \! \! \!
\end{equation}
\textcolor{black}{and therefore also a constant $c_a>0$ (not depending on $\epsilon$ or $\bar\delta_s$ or $\textcolor{black}{\bar \Delta_A}$ or $\textcolor{black}{\bar \Delta_B}$) such that
\begin{equation}\label{todev}\textstyle
    \mathcal D_\tau(\epsilon) \limsup\limits_{t\to +\infty}|\mathcal T(t,\zeta_t)|\le c_a\epsilon^\sharp
\end{equation}
for all $\epsilon\in (0,1)$, where $
\epsilon^\sharp$ was defined in \rref{epsilonsharp}, and where
\begin{equation}   \!  \! \!
    \begin{array}{l}\textstyle
    \mathcal D_\tau(\epsilon)=
1-\\\frac{2|K|\hat\tau\sqrt{2\epsilon}}{\sqrt{\pi}}  K_*|\textcolor{black}{\Delta_A}|_\infty \textcolor{black}{|B_0|_\infty}\!  \left(\lambda\sqrt{\frac{\xi_s}{\underline p}}\! +\! |B_{\rm new}|_\infty \sqrt{\frac{\epsilon}{2\pi}}\right).
    \end{array}\! \! \! \!  \! \!
\end{equation}
Using \rref{newepscond}, it follows that we can choose $c_{\mathcal T}=c_a/\mathcal D_\tau(\epsilon)$.} If we again use   \rref{xisbound1},  then we \textcolor{black}{obtain} the required constant $d_{\mathcal T}$, using \rref{kn10}.
Moreover,   we can substitute    \rref{MCD}-\rref{xis} to express our ultimate bound from \rref{kn18a} in terms of  $\epsilon$.\hfill$\square$

If  $A$ and $B$ in \rref{delayeq} are such that Assumptions \ref{as1}-\ref{as3} hold with $\bar\delta=0$ and $\bar\Delta_B=0$, then for any constant $\tau>0$, the product rule implies that Assumption \ref{as1} also holds with the same $A$ but  $B(t)$ replaced by  $B_{\rm new}(t)=\Phi_E(t-\tau,t)B(t+\tau)$. Therefore, in this case, Assumption \ref{as3}  holds for small enough $\epsilon>0$ \textcolor{black}{and $\bar\delta_s>0$} when we replace $B(t)$ by $B_{\rm new}(t)$ in Assumption \ref{as1}. Therefore, to ensure that we can apply Theorem \ref{prop1}, it suffices to find sufficient conditions under which the assumptions of Lemma \ref{sufficlem} hold with $B$ replaced by $B_{\rm new}$, because this will imply that the remaining  Assumption \ref{as2} also holds with $B$ replaced by $B_{\rm new}$. The following proposition provides the required sufficient conditions, using the   singular value notation   from Section \ref{intro}\textcolor{black}{, and where for simplicity, we assume that the delay $\tau$ is a constant (so $\Delta_\tau$ is the zero function in \rref{tdecomp}) and that  $\textcolor{black}{\Delta_B}$ is the zero function (but the general case follows from similar arguments)}:

 \begin{proposition}\label{as2remark}
 Let the assumptions of Lemma \ref{sufficlem} be satisfied with $A\in \mathbb R^{n\times n}$ being constant, and
for some   positive constants $\underline b$, $\bar A$, $\bar B$, $\bar \Gamma$, and $\Delta$ and
positive definite matrix $K\in \mathbb R^{n\times n}$. Let $\tau>0$ be \textcolor{black}{a given constant}, and set \begin{equation}\label{Lstar}\textstyle
     L_*=\left(\frac{\sigma_{\rm max}\left(e^{\scriptscriptstyle -\tau A}\right)}{\sigma_{\rm min}\left(e^{\scriptscriptstyle -\tau A}\right)}\right)^2.\end{equation} Assume that
\begin{equation}\label{sy18aaa}\asaa \begin{array}{l}
\textcolor{black}{|KA+A^\top K|_\infty}\left[\frac{1}{2} + \frac{2|K|^4\Delta^2\bar\Gamma^2\bar B^2}{
\underline{b}\lambda^2_{\rm min}(K)}L^3_*\right]
\\+
2 \bar{A} \Delta  \bar{\Gamma} |K|^2L_*- \frac{\underline{b}}{4}\lambda^2_{\rm min}(K) < 0.\end{array}
\end{equation}
Then the conclusions of Lemma \ref{sufficlem} are also satisfied with $B(t)$,
$\underline b$,  $\bar\Gamma$, and $K$ in Lemma \ref{sufficlem}
replaced by $B_{\rm new}(t)=e^{-\tau A}B(t+\tau)$, $\underline b_{\rm new}=\sigma^2_{\rm min}(e^{\scriptscriptstyle -\tau A})\underline b$,  $\bar\Gamma_{\rm new}=|e^{-\tau A}|^2\bar\Gamma$, and
$K_{\rm new}=K/\sigma^2_{\rm min}(e^{\scriptscriptstyle -\tau A})$, respectively.
\hfill$\square$\end{proposition}

{\em Proof.}
By assumption, our condition
\begin{equation}\label{sy18aa} \begin{array}{l}\textcolor{black}{|KA+A^\top K|_\infty}
\left[\frac{1}{2} + \frac{2|K|^4\Delta^2\bar\Gamma^2\bar B^2}{
\underline{b}\lambda^2_{\rm min}(K)}\right]
\\+
2 \bar{A} \Delta  \bar{\Gamma} |K|^2- \frac{\underline{b}}{4}\lambda^2_{\rm min}(K) < 0.\end{array}
\end{equation}from \rref{sy18} of Lemma \ref{sufficlem} is satisfied.
Also, since $A$ is constant, we have $\Phi_E(t-\tau,t)=e^{-\scriptscriptstyle \tau A}$. Therefore, our
 persistence of excitation condition \rref{PE} and   bound \rref{sy16} from Lemma \ref{sufficlem} hold with the new choice  $B_{\rm new}$(t) and with
 the new choices of $\underline b$ and $\bar\Gamma$ indicated above, using the fact that singular values are invariant under transposition, and the fact that the constantness of $A$ allows us to factor out  $e^{\scriptscriptstyle-\tau A}$ (resp.,  $e^{\scriptscriptstyle -\tau A^\top}$) from the left side (resp., right side) of the integrand in the persistence of excitation condition.
It follows that, with the   choice $B_{\rm new}$, the sufficient conditions from Lemma \ref{sufficlem} will hold if  the positive definite matrix $K_{\rm new}$   satisfies \begin{equation}\label{sy18a}\! \! \!   \!  \begin{array}{l}
\textcolor{black}{|K_{\rm new}A+A^\top K_{\rm new}|_\infty}\left[\frac{1}{2} + \frac{2|K_{\rm new}|^4\Delta^2\bar\Gamma^2\bar B^2 |e^{-\tau A}|^6}{
\underline{b}\sigma^2_{\rm min}(e^{-\tau A})\lambda^2_{\rm min}(K_{\rm new})}\right]
\\[.15em]+
2 \bar{A} \Delta  \bar{\Gamma}|e^{\scriptscriptstyle -\tau A}|^2 |K_{\rm new}|^2- \frac{\underline{b}}{4}\sigma^2_{\rm min}(e^{\scriptscriptstyle -\tau A})\lambda^2_{\rm min}(K_{\rm new})\\ < 0\end{array} \! \! \! \! \!
\end{equation}because $|B_{\rm new}|_\infty\le |e^{\scriptscriptstyle-\tau A}||B|_\infty\le |e^{\scriptscriptstyle-\tau A}|\bar B$.
However,
if we choose $K_{\rm new}=K/\sigma^2_{\rm min}(e^{\scriptscriptstyle -\tau A})$ in \rref{sy18a}, then we get
\begin{equation}\label{sy18a1} \asb\begin{array}{l}
\frac{\textcolor{black}{|KA+A^\top K|_\infty}}{
\sigma^2_{\rm min}(e^{-\tau A})}
\left[\frac{1}{2} + \frac{2|K|^4\Delta^2\bar\Gamma^2\bar B^2 |e^{-\tau A}|^6}{
\underline{b}\sigma^6_{\rm min}(e^{-\tau A})\lambda^2_{\rm min}(K)}\right]
\\+
\frac{2 \bar{A} \Delta  \bar{\Gamma}|e^{-\tau A}|^2 |K|^2}{\sigma^4_{\rm min}(e^{-\tau A})}- \frac{\underline{b}\lambda^2_{\rm min}(K)}{4\sigma^2_{\rm min}(e^{-\tau A})}
 < 0. \end{array}
\end{equation}
By multiplying   \rref{sy18a1}   by $\sigma^2_{\rm min}(e^{\scriptscriptstyle -\tau A})$,  we get  \rref{sy18aaa}, because the matrix norm equals the largest singular value. Hence, \rref{sy18a} is satisfied.
This concludes the proof. \hfill$\square$

Theorems \ref{thm1}-\ref{prop1}   allow cases where $A$ and $B$ are  unknown\mm{, provided we know bounds satisfying   Assumptions \ref{as1}-\ref{as2}}. On the other hand, our results are new, even in the basic scalar valued bounded extremum seeking case of $\dot x=u$, in which case   Assumptions \ref{as1}-\ref{as2} hold with $A(t)=0$ and $B(t)=B_{\rm new}(t)=1$ for all $t\ge 0$,
 $
 \bar A=0$, $\bar B=1$, $\bar D_B=0$, $P=K=1$, $\bar p=\underline p=1$, and $q=2$.
   With the preceding choices, we then have this special case of Theorem \ref{prop1} \textcolor{black}{with $\bar\Delta_A$, $\bar\Delta_B$, $\bar\Delta_\tau$, $\bar\delta$, and $\bar\delta_s=\bar{\mathcal M}_\Delta$ all being zero},  whose assumptions hold when $\epsilon>0$ is small enough for any choice of the constant delay $\tau>0$:

\begin{corollary}\label{cor1}
Let $\sigma_0>0$ and $\tau>0$ be any constants, and choose the constants
\begin{equation}\! \! \! \! \begin{array}{l}
c_1= 2\sqrt{2} \sqrt{\frac{\epsilon}{\pi}}\! \left(1\! +\! \frac{\epsilon }{4\pi}\right),\; \;  c_2=\frac{7\lambda\epsilon}{2\pi}
,\; \;  c_3=2\lambda^2\! \sqrt{\frac{2\epsilon}{\pi}} ,\\  \lambda=\frac{4\pi}{4\pi- \epsilon},\; \mm{
 d_1=
\frac{c_1}{2},\; d_2=\frac{c_1+c_3}{2}+c_2,\; d_3=\frac{c_3}{2},\\}
\xi_s = \frac{c_1}{
1+\sqrt{1-c_1c_3}},\; \,  \text{and}\; \,
\xi_l = \frac{1+\sqrt{1-c_1c_3}}{c_3},
    \end{array}\! \! \! \!
\end{equation}
    where $\epsilon\! \in\!  (0,2\pi  ]$ is any   constant such that
\begin{equation}\begin{array}{l}
c_1+c_3+2c_2\le 2,\;
 c_1c_3<1,\; \text{and}\\
2\left(1+\frac{\epsilon}{4\pi} \right)^2 \left(\sigma^2_0+\frac{\epsilon}{2\pi}\right)<\xi_l.
\end{array}\end{equation}
Then for each constant
\begin{equation}\label{xi011}\! \! \begin{array}{l}
    \xi_0\in \left(\max\left\{\xi_s,2 \left(1\! +\! \frac{\epsilon}{4\pi}  \right)^2\left(\sigma^2_0\! +\!  \frac{\epsilon}{2\pi} \right)\right\},\xi_l\right)\end{array}\! \! \!
\end{equation}
\textcolor{black}{and for  each initial state $x(0)$}   satisfying $|\zeta(0)|\le \sigma_0$,
the solution $x: [0,+\infty)\to \mathbb R$ of $\dot x(t)=u(t-\tau)$, in closed loop with the control
 \begin{equation}
\label{uchoice2aa}\textstyle
u(t) =
\sqrt{\frac{2\pi}{\epsilon}} \cos\left(\frac{2\pi}{\epsilon} t
+  \zeta^2(t)\right),
\end{equation}
where   $
\zeta(t) = x(t) +
\int_{t-\tau}^{t}
u(m)
{\rm d}m$ for all $t\ge 0$ and $\zeta(t)=\zeta(0)$ for all $t<0$,
satisfies
\begin{equation}\label{conclusion1d}\begin{array}{rcl} |x(t)|&\le& \lambda
    \sqrt{\frac{(\xi_0 - \xi_s) \xi_le^{-0.5t (\xi_l - \xi_s) c_3} + \xi_s (\xi_l - \xi_0) }{
(\xi_0 - \xi_s)e^{-0.5t (\xi_l - \xi_s) c_3}
+
\xi_l - \xi_0
}}\\&&+\,     \sqrt{\frac{\epsilon}{2\pi}}
+\left\vert \int_{t-\tau}^{t}
u(m)
{\rm d}m\right\vert\end{array}
\end{equation}
for all $t\ge 0$.
\hfill$\square$\end{corollary}

\begin{remark}
     Whereas standard reduction model controls allow arbitrarily long constant delays without any additional conditions on the delay, here we have  additional conditions on $\tau$  arising from  Assumption \ref{as2}, which has no analog in the standard theory of reduction model controls. Our additional  conditions on $\tau$ are that \rref{sy18aaa} must hold, which will hold when the corresponding undelayed system satisfies the sufficient conditions from Lemma \ref{sufficlem} and when $L_*$ is close enough to $1$, which covers cases with arbitrarily long delays $\tau>0$.
     We believe that our conditions on $\tau$ are  the price to pay  to achieve our feedback control goals for extremum seeking controls, which were beyond the scope of earlier reduction model approaches. In particular,  Proposition \ref{as2remark} covers scalar systems (where $L_*=1$)
and
 cases where the  singular values of $e^{\scriptscriptstyle -\tau A}$ are close enough together to satisfy   \rref{sy18aaa} when \rref{sy18aa} holds.
In Section \ref{illu}, we apply the preceding approach to a two dimensional system where $A$ is exponentially unstable.
\hfill$\square$\end{remark}

\section{Illustrations}
\label{illu}

%\textcolor{black}{This is how far MM revised  on 12/20.}

\subsection{Two Dimensional System}

We revisit the system
\begin{equation}
\label{f50}
\left[
\begin{array}{c}
\dot{x}_1
\\
\dot{x}_2
\end{array}
\right]
=
\left[
\begin{array}{cc}
2.1 & 4.9
\\
- 7.5 & 3.6
\end{array}
\right]
\left[
\begin{array}{c}
x_1
\\
x_2
\end{array}
\right]
+
\left[
\begin{array}{c}
\cos(\beta t)
\\
\sin(\beta t)
\end{array}
\right] u
\end{equation}
from the main example in \cite{SK13}, where  $\beta > 0$ is a constant. As noted in \cite{SK13}, a physical motivation for \rref{f50} is that $x=[x_1,x_2]^\top$ models the location of a mobile robot in the plane whose angular velocity actuator failed and is stuck at $\beta$, where the goal is stabilization to the origin using only the forward velocity $u$, in the presence of a position dependent perturbation term that is given by the drift term in \rref{f50}, i.e., $Ax$, which produces an unstable uncontrolled system having the poles $2.85\pm  10.7i$. This example is beyond the scope of works such as \cite{ZF24}, which treats the drift as a small perturbation term.

On the other hand, our Assumption \ref{as2} is satisfied with $K=\textcolor{black}{90}I$. To see why, note that with the choice
\begin{equation}
\label{f51}
B(t) = \left[
\begin{array}{c}
\cos(\beta t)
\\
\sin(\beta t)
\end{array}
\right],
\end{equation}we\mm{ get
\begin{equation}\! \!
\label{f52}\begin{array}{l}
B(t) B(t)^\top = \left[
\begin{array}{cc}
\cos^2(\beta t) & \sin(\beta t)\cos(\beta t)
\\
\sin(\beta t)\cos(\beta t) & \sin^2(\beta t)
\end{array}
\right] .\end{array}\!  \!
\end{equation}
We} \textcolor{black}{can   use the half angle  formulas for $\sin$ and $\cos$ to get}
\begin{equation}
\label{f54}
\frac{\beta}{2\pi} \int_t^{t + \frac{2\pi}{\beta}} B(m) B(m)^\top{\rm d}m = \frac{1}{2}\left[
\begin{array}{cc}
1 & 0
\\
0 & 1
\end{array}
\right] ,
\end{equation}
so we can choose $\Delta=2\pi/\beta$ and $\underline b=1/2$ to satisfy the PE requirement from \rref{PE} from  Lemma \ref{sufficlem}. Since we can also satisfy the requirement \rref{sy18} from Lemma \ref{sufficlem} for small enough $\Delta>0$, it follows that
  this example is covered by Theorem \ref{thm1}, when $\epsilon>0$, \textcolor{black}{$\bar\delta_s$,} and $\bar\delta$ are  small enough.

For instance, using the Mathematica computer program,
 we found that for $\sigma_0=1$, $\bar A=|A|$, $\bar B=1$, $\bar D_B=\beta$, \textcolor{black}{$K=90I$, $\bar\delta=\bar\delta_s=0$, and   $\beta= 8500$}, and using the weights \textcolor{black}{$a=0.5$ and $b=1.3$} in Remark \ref{weights},
 we can satisfy  our assumptions using $\epsilon=10^{\textcolor{black}{-6}}$. On the other hand, smaller choices of $\epsilon>0$ continue to satisfy our requirements from Assumption \ref{as3} and they produce smaller ultimate bounds on the right side of \rref{conc2}. We illustrate this point in \textcolor{black}{Table \ref{tab1a}, which contains} the results of our Mathematica calculations showing how the ultimate bound decreases as $\epsilon$ is chosen smaller, when $\sigma_0=1$ and all of the other parameter values   are chosen as above.
 \begin{table}[ht] \textcolor{black}{   \vspace{.2em}
    \hspace{4em}\scalebox{1}{\begin{tabular}{|l|c|c|}\hline
      {\bf  $\epsilon$ in Control \rref{uchoice}}  & $10^{-6}$ &$10^{-7}$
      \\\hline
      {\bf Ultimate Bound}& 0.291879 & 0.16251
      \\\hline
    \end{tabular}}} \vspace{.25em}
    \caption{Relationship Between $\epsilon$ and Ultimate Bound for \rref{f50} in Theorem \ref{thm1} with $K=\textcolor{black}{90}I$ } \label{tab1a} \vspace{-.5em}
\end{table}

\mm{\begin{table}[ht]
    \centering
    \centering\scalebox{1}{\begin{tabular}{|l|c|c|}\hline
      {\bf Choice of $\epsilon$ in Control \rref{uchoice}}   & $10^{-9}$& $10^{-10}$
      \\\hline
      {\bf Ultimate Bound from \rref{conc2}}&0.0149648& 0.00841221
      \\\hline
    \end{tabular}} \vspace{.25em}
    \caption{Relationship Between $\epsilon$ and Ultimate Bound for \rref{f50} in Theorem \ref{thm1} with $K=300I$ for Smaller $\epsilon$ Values } \label{tab1b} \vspace{-.05em}
\end{table}}

This is further illustrated in Fig. 1, which shows Mathematica plots for  $x_1$ and $x_2$ for  \rref{a2v} with the controller from Theorem \ref{thm1} and the previous parameter values \textcolor{black}{with $\epsilon=10^{-6}$, for three  initial states}.\mm{ As illustrated in the figures, the \textcolor{black}{smaller value of $\epsilon$    produced the faster  and closer} convergence to the zero equilibrium.}
\begin{figure}[!htb]
\vspace{-.07em}\centering
  \includegraphics[scale=0.417]{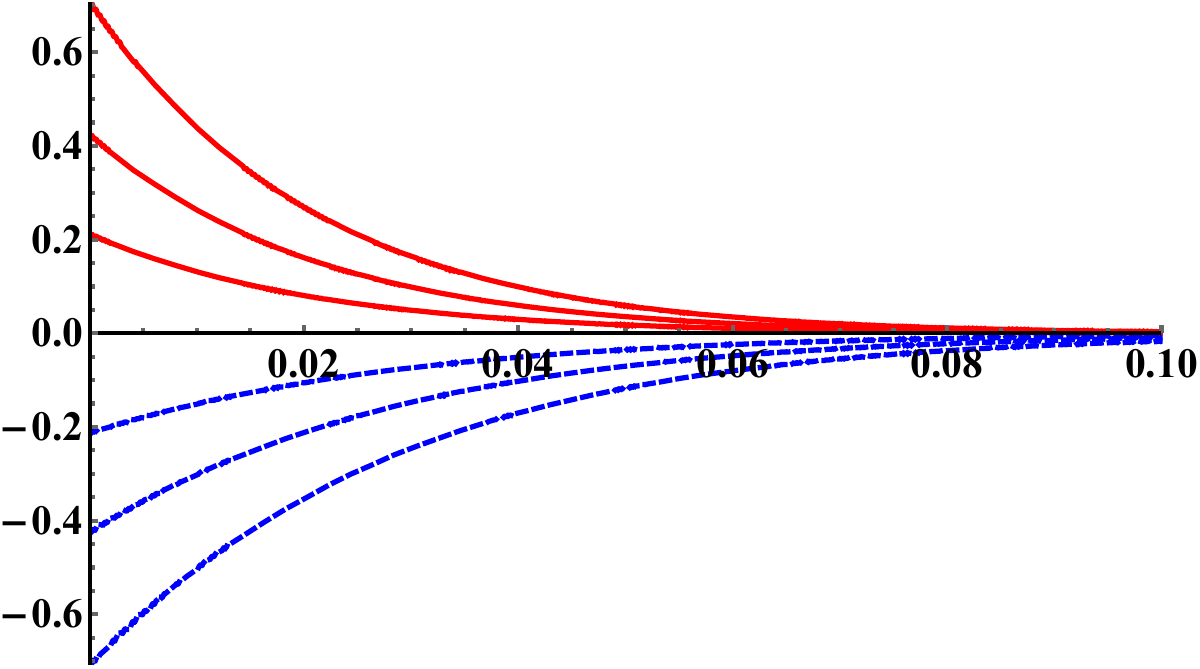}  \vspace{-1.95em}\caption{Closed Loop Solutions for \rref{f50} with Initial \textcolor{black}{States $x(0)=(1/\sqrt{2},-1/\sqrt{2})$,
  $x(0)=(0.6/\sqrt{2},-0.6/\sqrt{2})$, and $x(0)=(0.3/\sqrt{2},-0.3/\sqrt{2})$
  Showing $x_1(t)$ (Solid Red) and $x_2(t)$ (Dashed Blue)  on $0.1$} Second Time Horizon}\vspace{-.3em} \label{firstfig}\end{figure}

If  we replace the weight $\textcolor{black}{a=0.5}$ from Remark \ref{weights} by $a=1$, and keep all other parameters the same in this example, then instead of the ultimate bounds in the second \textcolor{black}{line of Table \ref{tab1a}, we  obtained $0.417193$ and  $0.230003$} for the indicated values of $\epsilon$ in the first line of Table \ref{tab1a}, respectively. This illustrates how a suitable choice of the weight $a>0$ can reduce the ultimate bounds by more than one \textcolor{black}{third}.

When $\bar\Delta_A=\bar\Delta_B=0$, we can also apply our reduction model result from Theorem   \ref{prop1} to \rref{f50}, with the special case
    $B_{\rm new}(t)=e^{-\tau A}B(t+\tau)$ of $B_{\rm new}$ from \rref{94} for cases where $A$  \textcolor{black}{and $\tau$ are constant (so $\Delta_\tau$ is the zero function in \rref{tdecomp})}.
    We illustrate this point in Table \ref{tab1d}, which shows  the ultimate bound on $\zeta(t)$  \textcolor{black}{for  $\tau=0.05$}, for two choices of the $\textcolor{black}{\epsilon}>0$, where all of the other parameters were kept the same as in the undelayed case, except with \textcolor{black}{
    $K=70I$} in Proposition \ref{as2remark} and with
    the weights \textcolor{black}{$a=0.5$ and $b=1.06$.\mm{ This compares favorably with the maximum delay bound of $\tau_M=0.0050$ from \cite[Table 1]{ZF24} that was obtained using the previous averaging approaches.}}
\begin{table}[ht] \vspace{.75em}\textcolor{black}{
    \hspace{4em}\scalebox{1}{\begin{tabular}{|l|c|c|}\hline
      {\bf   $\epsilon$ in Control \rref{uchoice}}  & $10^{-6}$ &$10^{-7}$
      \\\hline
      {\bf Ultimate Bound}& 0.32505 & 0.178749
      \\\hline
    \end{tabular}}} \vspace{.25em}
    \caption{Relationship Between $\epsilon$ and Ultimate Bound for \rref{f50} \textcolor{black}{in Theorem \ref{prop1} with $K=\textcolor{black}{70}I$ \textcolor{black}{and $\tau=0.05$}}} \label{tab1d}
\end{table}
\textcolor{black}{However,} by  reducing $|K|$ \textcolor{black}{or $\epsilon$}, we can allow larger constant $\tau>0$ values.

For instance, with $\tau=1$, our assumptions are satisfied for this example with $K=\textcolor{black}{36}I$ in Proposition \ref{as2remark}, \textcolor{black}{$\beta=10000$,} $\epsilon=10^{-8}$, \textcolor{black}{$a=0.1$, $b=1$},  and all other parameters as we selected them in Table \ref{tab1d}, and these choices produced the ultimate bound \textcolor{black}{$0.0169173$}. This compares favorably with the delay bound that we would have had, if we instead used the method from Section \ref{delays}, which would have called for choosing $\bar\delta\le 0.001$ when $\epsilon=10^{-7}$ (and using the other parameter values we used for Table \ref{tab1a}) and then applying \rref{taubound} to obtain the delay bound (which would have required $|\tau|_\infty\le 10^{-10}$).

\textcolor{black}{Finally, we can cover nonzero added uncertainties $\Delta_\tau$, since for instance, the assumptions of Theorem \ref{prop1} hold for the preceding example with $\beta=10000$, $\hat\tau=0.5$, \textcolor{black}{$\bar\Delta_A=\bar\Delta_B=0$,} $k=100$, $\epsilon=10^{-8}$, and $|\Delta_\tau|_\infty\le 6\times 10^{-4}$. This bound  is consistent with the value of $\epsilon$.}
This illustrates how our methods   apply for reduction model  controllers to quantify the effects of time delays that were beyond the scope of earlier extremum seeking  results.

\mm{This example \rref{f50} also satisfies our assumptions using $K$ matrices that are not scalar matrices, and we found that different choices of $K$ can reduce the ultimate bound, as compared with the values that we obtained in Table \ref{tab1}. For instance, if we change the scalar matrix $K$ that was used to generate Table \ref{tab1} to the nonscalar matrix
\begin{equation}\label{newS}
    K=100\left[\begin{array}{cc}1&0\\0.5&1\end{array}\right]
\end{equation}
and keep all other values the same as the ones that we used to generate Table \ref{tab1}, then we obtain the smaller ultimate bounds for each of the four choices of $\epsilon$ that we show in Table \ref{tab1a}. This suggests the problem of optimizing over $K$ to obtain the largest possible ultimate bounds under the assumptions of this paper, which we believe is an open problem.
\begin{table}[ht]
    \centering
    \centering\scalebox{1}{\begin{tabular}{|l|c|c|c|c|}\hline
      {\bf Choice of $\epsilon$ in Control \rref{uchoice}}  & $10^{-7}$ &$10^{-8}$& $10^{-9}$& $10^{-10}$
      \\\hline
      {\bf Ultimate Bound from \rref{conc2}}& 0.0670852 & 0.037663&0.0211679& 0.0119004
      \\\hline
    \end{tabular}}\vspace{-.5em}
    \caption{Relationship Between $\epsilon$ and Ultimate Bound for \rref{f50} in Theorem \ref{thm1} with \rref{newS} \label{tab1a}}  \vspace{-.75em}
\end{table}}

\subsection{One Dimensional System}

Given a bounded continuous function $\mathcal A:[0,+\infty)\to \mathbb R$,
we consider the system
\begin{equation}
\label{h}
\dot{x}(t) =\mathcal A(t) x + \cos(\beta t) u
\end{equation}
with a constant $\beta > 0$. This includes the one dimensional example in \cite{SK13}, in the special case where $\mathcal A(t)=1$ for all $t\ge 0$.
Then, with $K = k$ for any positive constant $k$, our choice of $H$  \textcolor{black}{in \rref{a17v}} gives
\begin{equation}
\label{h1}\textstyle
H(t) = \mathcal A(t) - k \cos^2(\beta t) = \mathcal A(t) - \frac{k}{2} - \frac{k}{2}\cos(2 \beta t).
\end{equation}
One can then use Lemma \ref{sufficlem} to find the required   function $P(t)$, because in this case, the persistence of excitation condition \rref{PE} from the lemma has the form
\begin{equation}\label{PEs}\textstyle
\frac{1}{\Delta}\int_t^{t+\Delta} \cos^2(\beta \ell){\rm d}\ell\ge \underline b,
    \end{equation}
which holds with $\Delta=2\pi/\beta$ and $\underline b=1/2$.

However, we find it more convenient  and less restrictive to find $P$ in the following alternative way (but see Remark \ref{floq} below, where $P$ is found using   Lemma \ref{sufficlem}).
Let
\begin{equation}
\label{h2}
P(t) = e^{ \frac{k}{2\beta} \sin(2\beta t)}.
\end{equation}
Then
\begin{equation}
\label{h3}\! \! \! \ashhneg \begin{array}{rcl}
\dot{P}(t) + 2 H(t) P(t) &=&  k e^{ \frac{k}{2\beta} \sin(2\beta t)} \cos(2\beta t)
\\&&+ 2 \left[\mathcal A(t) - \frac{k}{2} - \frac{k}{2}\cos(2 \beta t)\right]\\&&\; \; \; \times e^{ \frac{k}{2\beta} \sin(2\beta t)}\\
&
=& \left(2\mathcal A(t) - k\right) P(t).\end{array}\! \! \!
\end{equation}
Then Assumption \ref{as2} is satisfied for any $\beta > 0$ and any $k>2|\mathcal A|_\infty$, since we also have
\begin{equation}
\label{a}
 e^{- \frac{k}{2\beta}} \leq P(t) \leq e^{\frac{k}{2\beta}}
\; \text{and}\;
|\dot{P}(t)| \leq k e^{\frac{k}{2\beta}}
\end{equation}for all $t\ge 0$.
This allows us to apply Theorem \ref{thm1}. For instance, if we choose     $\sigma_0=1$, $\mathcal A(t)=1$ for all $t\ge 0$, \textcolor{black}{$\bar\delta=\bar\delta_s=0$}, $\beta=500$,  $K=20$, and  the weights in Remark \ref{weights} to be $a=0.1$ and $b=1$, then we can satisfy our assumptions with $\epsilon=10^{-5}$. Reducing $\epsilon>0$ while keeping all other parameters the same can reduce the ultimate bound. In Tables \ref{tab3a}-\ref{tab3b}, we illustrate the effects of reducing $\epsilon$ in this example, with the preceding parameter values.

We can also allow  $\sigma_0>1$. For instance, with $\epsilon=10^{-8}$ but all other parameters   the same as in Table \ref{tab3a}, we used   Mathematica to find that the largest allowable $\sigma_0$ for which our assumptions are all satisfied was $\sigma_0=10.87$, which gave the same ultimate bound $0.0136521$.
\begin{table}[ht]
    \centering
    \centering\scalebox{1}{\begin{tabular}{|l|c|c|c|c|}\hline
      {\bf Choice of $\epsilon$ in Control \rref{uchoice}}  & $10^{-5}$ &$10^{-6}$
      \\\hline
      {\bf Ultimate Bound from \rref{conc2}}& 0.090951 & 0.0503611
      \\\hline
    \end{tabular}}\vspace{.25em}
    \caption{Relationship Between $\epsilon$ and Ultimate Bound in One-Dimensional Example in \cite{SK13} }\label{tab3a}
     \vspace{-1.25em}
\end{table}

\begin{table}[ht]
    \centering\vspace{.25em}
    \centering\scalebox{1}{\begin{tabular}{|l|c|c|}\hline
      {\bf Choice of $\epsilon$ in Control \rref{uchoice}}  & $10^{-7}$& $10^{-8}$
      \\\hline
      {\bf Ultimate Bound from \rref{conc2}} &0.0282206& 0.0158385
      \\\hline
    \end{tabular}}\vspace{.25em}
    \caption{Relationship Between $\epsilon$ and Ultimate Bound in One-Dimensional Example in \cite{SK13} for Smaller $\epsilon$ Values }\label{tab3b}\vspace{-.25em}
\end{table}

If we  use   $a=1$ instead of $a=0.1$ in Remark \ref{weights} and keep all other parameters the same with $\sigma_0=1$, then the ultimate bounds from \rref{conc2} are instead $0.283668$, $0.158466$, $0.0889727$, and $0.0499997$ for the corresponding $\epsilon$ values from the first lines of Tables \ref{tab3a}-\ref{tab3b}, respectively. This   illustrates how choosing $a$ differently from $1$ can reduce the ultimate bounds by more than one half.
We can also satisfy our assumptions by reducing $\sigma_0$ to values in $(0,1)$, but this did not change the largest allowable value $10^{-5}$ for  $\epsilon$   for which our assumptions hold  when all other parameter values were kept the same.
\begin{remark}\label{floq}
Our choice \rref{h2} of   $P$  is a special case of Floquet theory constructions. An alternative  approach    would have been to solve the separable differential equation
\begin{equation}\dot P+2P(\mathcal A(t)-k\cos^2(\beta t))=-qP\end{equation} for the dependent variable $P$ and for suitable positive constants $q$, which in the special case where $\mathcal A(t)=1$ for all $t\ge 0$  would have instead produced the general solution
\begin{equation}
P(t)=ce^{-t(q+2-k)}e^{\frac{k}{2\beta}\sin(2\beta t)}
\end{equation}for real constants $c$,
which would have violated  Assumption \ref{as2} unless $q=k-2$. We could also have constructed $P$   using Lemma \ref{sufficlem}, whose persistence of excitation condition \rref{PEs}
and other conditions hold in this example with $\underline b=\bar\Gamma=1/2$, $\Delta=2\pi/\beta$, and any   $k>0$ such that
\begin{equation}\label{kcond}\textstyle
    \frac{1}{k}+\frac{8\pi^2k}{\beta^2}+\frac{2\pi}{\beta}<\frac{1}{8}.
\end{equation}
Condition \rref{kcond} places more restrictions
on  $k$,  compared with
the approach we used above. \textcolor{black}{This is because the method that we used above uses the requirement $k>2|\mathcal A|_\infty$, but \rref{kcond} may be violated for large enough $k$ or small enough $\beta$. For instance, if we use the same parameter values that we used in Tables \ref{tab4a}-\ref{tab4b} except we reduce $\beta$ and $k$ to $\beta=200$ and $k=15$ respectively, then the assumptions of Theorem \ref{thm1} continue to be satisfied (with the ultimate bounds $0.0707977$, $0.0394922$, $0.0221094$, and $0.0124019$ for $\epsilon=10^{-5}$, $\epsilon=10^{-6}$, $\epsilon=10^{-7}$, and $\epsilon^{-8}$, respectively, instead of the ones in Tables \ref{tab4a}-\ref{tab4b}) but \rref{kcond} would be violated (because its left side would be $0.127691>1/8$). This motivates using Assumption \ref{as2} (instead of using the more restrictive assumptions from Lemma \ref{sufficlem}) in Theorem \ref{thm1}.}
\hfill$\square$\end{remark}

\subsection{Illustration of Corollary \ref{cor1}}
\textcolor{black}{When $\Delta_A$, $\Delta_B$, and $\Delta_\tau$ are   the zero function, the bound \rref{fTb} from o}ur proof of Theorem \ref{prop1} provides a constant $c_*>0$ (not depending on $\tau>0$ or on $\epsilon\in (0,1)$) such that $\limsup_{t\to +\infty}|\mathcal T(t,\xi_t)|\le c_* \epsilon^{\scriptscriptstyle 1/4}(\tau+1)$ for all choices of $\epsilon\in (0,1)$ and $\tau>0$. It follows from our bound on $|x(t)|$ from \rref{kn10}  that larger values of   $\tau>0$ can call for smaller values of $\epsilon$,   to achieve smaller ultimate bounds on the norm $|x(t)|$ of the state of the closed loop system. We next illustrate this point in the
 special case of $\dot x=u$ from Corollary \ref{cor1}.

 In  Tables \ref{tab4a}-\ref{tab4b} below, we show Mathematica calculations of
the ultimate bound for $|x(t)$ for different choices of the delay $\tau>0$ and $\epsilon>0$, in the special case where  $\sigma_0=1$ and $P=K=1$, and using $a=0.1$ and $b=0.9$ to obtain the weights in the formulas from Remark \ref{weights}.
\begin{table}[ht]   \vspace{.35em}
    \centering
    \centering\scalebox{1}{\begin{tabular}{|l|c|c|c|c|}\hline
      {\bf Choice of $\epsilon$ in Control}  & $10^{-2}$ &$10^{-3}$
      \\\hline
      {\bf Delay $\tau$ }& $10^{-2}$ & 0.02
      \\\hline
      {\bf Ultimate Bound}& 0.29442 & 0.118344\\\hline
    \end{tabular}}\vspace{.25em}
    \caption{Relationship Between Allowable $\epsilon$ and $\tau$ and Ultimate Bound for $\dot x=u$ }\label{tab4a}\vspace{-.95em}
\end{table}

\begin{table}[ht]  \vspace{.73em}
    \centering
    \centering\scalebox{1}{\begin{tabular}{|l|c|c|c|c|}\hline
      {\bf Choice of $\epsilon$ in Control}  & $10^{-4}$ &$10^{-5}$
      \\\hline
      {\bf Delay $\tau$ }& 0.03 & 0.04
      \\\hline
      {\bf Ultimate Bound}& 0.0520498 & 0.0249277\\\hline
    \end{tabular}}\vspace{.25em}
    \caption{Relationship Between Allowable $\epsilon$ and $\tau$ and Ultimate Bound for $\dot x=u$ for Smaller $\epsilon$ Values }\label{tab4b}
     \vspace{.35em}
\end{table}

\mm{As in the previous examples from \cite{SK13} from above, t}\textcolor{black}{This} illustrates how the ultimate bound decreases as we reduce $\epsilon>0$. This example is beyond the scope of \cite{MMF24a}, which was confined to two-dimensional systems under suitable conditions on the different weights in the control components.
Larger $\epsilon$ values can be selected in the simulations, to find upper bounds on the $\epsilon$ for which practical stability is preserved to address any conservativeness in the $\epsilon$ conditions in our theorems.

\section{Conclusions}
We provided a new mathematical analysis for  bounded extremum seeking  for a large class of linear time-varying systems with measurement uncertainty. Significantly novel features   included our new comparison function argument, which led to practical exponential stability bounds that quantify the effects of the small parameter. This \textcolor{black}{covered} unknown control directions, including cases where both coefficient matrices are subjected to uncertainty and   the uncontrolled systems can be unstable. We applied our method to construct reduction model controllers, which made it possible to quantify the effects of input delays.\mm{ Our examples illustrated benefits of our method\mm{ as compared with prior methods}.} Since unknown control directions frequently arise in aerial applications, we aim to apply our methods to aerial models that contain measurement uncertainty\textcolor{black}{, and to find least conservative upper bounds on the   parameters and least conservative ultimate bounds. We also aim to find   analogs for Newton-based extremum seeking with uncertainty estimation, which would provide analogs of works such as \cite{G20a} that would apply for bounded extremum seeking with uncertain delays and uncertain coefficient matrices}.

\bibliographystyle{abbrv}
\bibliography{Aug2024ESRefs}

\end{document}